\documentclass[12pt]{article}

\usepackage{amsfonts,amsthm,amsmath,latexsym,bm,graphics,indentfirst,amssymb,color}
\usepackage[hypertex]{hyperref}

\def\theequation{\thesection.\@arabic\c@equation}
\makeatother

\renewcommand{\theequation}{\thesection.\arabic{equation}}
\newtheorem{lemma}{Lemma}[section]

\newtheorem{proposition}{Proposition}[section]
\newtheorem{corollary}{Corollary}[section]
\newtheorem{remark}{Remark}[section]

\newtheorem{theorem}{Theorem}[section]
\newcommand{\ve} {\varepsilon}
\newcommand{\R}{{\mathbb R}}

\def\bs{\begin{split}}
\def\es{\end{split}}
\numberwithin{equation}{section}
\newcommand{\vect}[2]{\ensuremath{\left( \begin{array}{c}
            #1 \\
            #2
            \end{array}
        \right)}}

\title{On Non-topological  Solutions of the ${\bf G}_2$ Chern-Simons System}
\author{Weiwei Ao \footnote{Center for Advanced Study in Theoretical Science, National Taiwan University, Taipei, Taiwan.
{\sl weiweiao@gmail.com}} \; Chang-Shou Lin \footnote{Taida Institute of Mathematics, Center for Advanced study in Theoretical Science, National Taiwan University, Taipei, Taiwan. {\sl cslin@math.ntu.edu.tw}}\; Juncheng Wei\footnote{Department of Mathematics, University of British Columbia, Vancouver, BC V6T 1Z2 and Department of Mathematics, Chinese University of Hong Kong, Shatin, Hong Kong. {\sl jcwei@math.ubc.ca}}}

\begin{document}
\maketitle
\begin{abstract}
For any rank 2 of simple Lie algebra,
the relativistic Chern-Simons system has the following form:
\begin{equation}\label{e001}
\left\{\begin{array}{c}
\Delta u_1+(\sum_{i=1}^2K_{1i}e^{u_i}
-\sum_{i=1}^2\sum_{j=1}^2e^{u_i}K_{1i}e^{u_j}K_{ij})=4\pi\displaystyle \sum_{j=1}^{N_1}\delta_{p_j}\\
\Delta u_2+ (\sum_{i=1}^2K_{2i}e^{u_i}-\sum_{i=1}^2\sum_{j=1}^2e^{u_i}K_{2i}e^{u_j}K_{ij})=4\pi\displaystyle \sum_{j=1}^{N_2}\delta_{q_j}
\end{array}
\right.\mbox{in}\; \mathbb{R}^2,
\end{equation}
where $K$ is the Cartan matrix of rank $2$. There are three Cartan matrix of rank 2: ${\bf A}_2$, ${\bf B}_2$ and ${\bf G}_2$.
A long-standing open problem for \eqref{e001} is the question of the existence of non-topological solutions. In a previous paper \cite{ALW}, we have proven the existence of non-topological solutions for the ${\bf A}_2$ and ${\bf B}_2$ Chern-Simons system. In this paper, we continue to consider the ${\bf G}_2$ case. We prove the existence of non-topological solutions under the condition that either
$N_2\displaystyle\sum_{j=1}^{N_1} p_j=N_1\displaystyle \sum_{j=1}^{N_2} q_j $ or
$N_2\displaystyle\sum_{j=1}^{N_1}p_j \not =N_1\displaystyle \sum_{j=1}^{N_2} q_j$ and  $N_1,N_2>1$, $
 |N_1-N_2|\neq 1$. We solve this problem by a perturbation from the corresponding ${\bf G}_2$ Toda system with one singular source. Combining with \cite{ALW}, we have proved the existence of non-topological solutions to the Chern-Simons system with Cartan matrix of rank $2$.
\end{abstract}
\section{Introduction}
\subsection{Background}
There are four types of simple non-exceptional Lie Algebra:${\bf A}_m$, ${\bf B}_m$, ${\bf C}_m$, and ${\bf D}_m$ which Cartan subalgebra are $sl(m+1)$, $so(2m+1)$, $sp(m)$, and $so(2m)$ respectively.
{To each of them, a Toda system is associated}. In geometry, solutions of Toda system is closely related to holomorphic curves in projective spaces.
For example, the Toda system of type ${\bf A}_m$ can be derived from the classical Pl\"{u}cker formulas, and any holomorphic curve gives rise to a solution $u$ of the Toda system,
whose branch points correspond to the singularities of $u$. Conversely, we could integrate the Toda system, and any solution $u$ gives rise to a holomorphic curve in $\mathbb{CP}^n$ at least locally.
See \cite{da}, \cite{lwy} and reference therein. It is very interesting to note that the reverse process holds globally if the domain for the equation is $\mathbb{S}^2$ or $\mathbb{C}$. Any solution $u$ of type ${\bf A}_m$ Toda system
on $\mathbb{S}^2$ or $\mathbb{C}$ could produce a global holomorphic curves into $\mathbb{CP}^n$. This holds even when the solution $u$ has singularities. We refer the readers to \cite{lwy} for more precise statements of these results.

\medskip

In physics, the Toda system also plays an important role in non-Abelian gauge field theory. {One example} is the relativistic Chern-Simons model proposed by Dunne \cite{d0,d,d2} in order to explain the physics of high critical temperature superconductivity. See also \cite{kl}, \cite{l} and \cite{l2}.

\medskip

The model is defined in the (2+1) Minkowski space $\mathbb{R}^{1,2}$, the gauge group is a compact Lie group with a semi-simple Lie algebra $\mathcal{G}$. The Chern-Simons Lagrangian density $\mathcal{L}$ is defined by:
\begin{equation*}
\mathcal{L}=-k\epsilon^{\mu\nu\rho}tr(\partial_\mu A_\nu A_\rho+\frac{2}{3}A_\mu A_\nu A_\rho)-tr((D_\mu \phi)^\dag D^\mu \phi)-V(\phi,\phi^\dag)
\end{equation*}
for a Higgs field $\phi$ in the adjoint representation of the compact gauge group $G$, where the associated semi-simple Lie algebra is denoted by $\mathcal{G}$ and the $\mathcal{G}-$valued gauge field $A_\alpha$ on $2+1$ dimensional Minkowski space $\R^{1,2}$ with metric diag\{-1,1,1\}. Here $k>0$ is the Chern-Simons coupling parameter, tr is the trace in the matrix representation of $\mathcal{G}$ and $V$ is the potential energy density of the Higgs field $V(\phi,\phi^\dag)$ given by
\begin{equation*}
V(\phi,\phi^\dag)=\frac{1}{4k^2}tr(([[\phi,\phi^\dag],\phi]-v^2\phi)^\dag([[\phi,\phi^\dag],\phi]-v^2\phi)),
\end{equation*}
where $v>0$ is a constant which measures either the scale of the broken symmetry or the subcritical temperature of the system.

\medskip

In general, the Euler-Lagrangian equation corresponding $\mathcal{L}$ is very difficult to study. So we restrict to consider solutions to be energy minimizers of the Lagrangian functional, and then a self-dual system of first order derivatives could be derived from minimizing the energy functional:
\begin{equation}\label{selfdual}
\begin{array}{l}
D_{-}\phi=0\\
F_{+-}=\frac{1}{k^2}\left[v^2\phi-[[\phi , \phi ^{\dagger}],\phi ],\phi ^{\dagger}\right],
\end{array}
\end{equation}
where $D_{-}=D_1-iD_2$, and $F_{+-}=\partial _+A_--\partial _-A_++[A_+,A_-]$ with $A_{\pm}=A_1\pm iA_2$, $\partial _{\pm}=\partial_1 \pm i\partial _2$. Here $\partial _i$ and $D_i$ are respectively the partial derivative and the gauge-covariant derivative w.r.t $z_i$, $i=1,2$.

\medskip

In order to find non-trivial solutions which are not algebraic solutions of $[[\phi ,\phi ^{\dagger}],\phi ]=v^2\phi $, Dunne \cite{d} has considered a simplified form of the self-dual system \eqref{selfdual} in which both the gauge potential $A$ and the Higgs field $\phi $ are algebraically restricted, for example, $\phi $ has the following form:
\begin{align*}
\phi=\displaystyle \sum _{a=1}^r\phi ^aE_a,
\end{align*}
where $r$ is the rank of the Lie algebra $\mathcal{G}$, $\{E_{\pm a}\}$ is the family of the simple root step operators (with $E_{-a}=E_{a}^+$), and $\phi ^a$ are complex-valued functions.

\medskip

In this paper, we consider the case of rank 2. Let
\begin{align*}
u_a=\ln |\phi ^a|^2.
\end{align*}
By using this ansatz, then equation \eqref{selfdual} can be reduced to
\begin{equation}\label{equation1}
\left\{\begin{array}{c}
\Delta u_1+\frac{v^4}{k^2}(\sum_{i=1}^2K_{1i}e^{u_i}-\sum_{i=1}^2\sum_{j=1}^2e^{u_i}K_{1i}e^{u_j}K_{ij})=4\pi\displaystyle \sum_{j=1}^{N_1}\delta_{p_j}\\
\Delta u_2+\frac{v^4}{k^2} (\sum_{i=1}^2K_{2i}e^{u_i}-\sum_{i=1}^2\sum_{j=1}^2e^{u_i}K_{2i}e^{u_j}K_{ij})=4\pi\displaystyle \sum_{j=1}^{N_2}\delta_{q_j}
\end{array}
\right.\mbox{in}\; \mathbb{R}^2,
\end{equation}
where $K=\left(\begin{array}{cc}
K_{11}&K_{12}\\
K_{21}&K_{22}
\end{array}
\right)$ is the Cartan matrix of rank 2 of the Lie algebra $\mathcal{G}$, $\{p_1,\ldots, p_{N_1}\}$ and $\{ q_1,\ldots,q_{N_2}\}$ are given vortex points. For the details of the process to derive \eqref{equation1} from (\ref{selfdual}), we refer to \cite{d},\cite{nt2},\cite{y} and \cite{y2}. In this paper, without loss of generality, we assume $\frac{v^4}{k^2}=1$.

\medskip

It is known that there are only three types of Cartan matrix  of rank 2, given by
\begin{equation}
\label{cartan1}
{\bf A}_2= \begin{pmatrix} 2 & -1 \\
-1 & 2
\end{pmatrix} \ , \ {\bf B}_2(={\bf C}_2)= \begin{pmatrix} 2 & -1 \\
-2 & 2
\end{pmatrix} \ , \ {\bf G}_2= \begin{pmatrix} 2 & -1 \\
-3 & 2
\end{pmatrix} \  .
\end{equation}

\medskip

In the previous paper \cite{ALW}, we have constructed non-topological solutions in the case of ${\bf A}_2$ and ${\bf B}_2$. In this paper, we will construct non-topological solutions for the ${\bf G}_2$ case, i.e. the following equation:
\begin{equation}\label{eqofg}
\left\{\begin{array}{c}
\Delta u_1+2e^{u_1}-e^{u_2}=4e^{2u_1}-2e^{2u_2}+e^{u_1+u_2}+4\pi\displaystyle\sum_{j=1}^{N_1}\delta_{p_j}\\
\Delta u_2+2e^{u_2}-3e^{u_1}=4e^{2u_2}-6e^{2u_1}-3e^{u_1+u_2}+4\pi\displaystyle\sum_{j=1}^{N_2}\delta_{q_j}.
\end{array}
\right.
\end{equation}

\subsection{Previous Results}

In the literature, a solution ${\bf u}=(u_1,u_2)$ to system (\ref{equation1}) is called a {\it topological} solution if ${\bf u}$ satisfies
\begin{align*}
{u_a(z)\rightarrow \ln\sum_{j=1}^2(K^{-1})_{aj} }\quad \text{as}\; |z|\rightarrow +\infty, \;a=1,2,
\end{align*}
and is called a {\it non-topological} solution if ${\bf u}$ satisfies
\begin{align}
u_a(z)\rightarrow -\infty \quad \text{as}\; |z|\rightarrow +\infty, \;a=1,2.
\end{align}

\medskip

The existence of topological solutions with arbitrary multiple vortex points was proved by Yang \cite{y2} more than fifteen years ago, not only for Cartan matrix of rank $2$, but also for general Cartan matrix including $SU(N+1)$ case, $N\geq 1$. However, the existence of non-topological solutions is more difficult to prove.   The first result was due to
 Chae and Imanuvilov \cite{ci} for the $SU(2)$ Abelian Chern-Simons equation which is obtained by letting $u_1(z)=u_2(z)=u(z)$ in the ${\bf A}_2$ system  where $u$ satisfies
\begin{equation}\label{e103}
\Delta u+e^u(1-e^u)=4\pi \displaystyle \sum _{j=1}^N \delta_{p_j}\mbox{ in }\; \mathbb{R}^2.
\end{equation}

\medskip

Equation \eqref{e103} is the $SU(2)$ Chern-Simons equation for the Abelian case. This relativistic Chern-Simons model was proposed by Jackiw-Weinberg \cite{jw} and Hong-Kim-Pac \cite{hkp}. For the past more than twenty years, the existence and multiplicity of solutions to \eqref{e103} with different nature (e.g. topological, non-topological, periodically constrained etc.) have been studied, see \cite{cy}, \cite{ci}, \cite{cfl}, \cite{ccl}, \cite{ckl}, \cite{hkp}, \cite{ly}, \cite{ly2}, \cite{mr}, \cite{mr1}, \cite{sy1}, \cite{sy2}, \cite{t}, \cite{t2}, \cite{w} and references therein.

\medskip

In \cite{ci},  Chae and Imanuvilov
proved the existence of non-topological solutions for (\ref{e103}) for any  vortex points $(p_1, ..., p_N)$. {For the question of existence of non-topological solutions for the ${\bf A}_2$ system , an {``answer"}  was given by Wang and Zhang \cite{wz} but their proof contains serious gaps.  In fact   they used a special  solution of the Toda system as the approximate solution, but they did not have the full non-degeneracy of the linearized equation of the Toda system and  their analysis for the linearized equation is incorrect.} Thus, the existence of non-topological solutions has remained a long-standing open problem. Even for radially symmetric solutions (the case when all the vortices coincide), the ODE system  is much more subtle than equation \eqref{e103}. The classification of radial solution is an important issue for future study as long as bubbling solutions are concerned, see \cite{ci}, \cite{cfl}, \cite{ccl}, \cite{ckl}, \cite{edm}, \cite{ly}, \cite{ly2}, \cite{nt}, \cite{pt},  \cite{t}, \cite{t2}, \cite{wzz} in this direction.

\medskip

For the rank $2$ Chern-Simons system of Lie type, Huang and the second author \cite{HL, HL1} studied the structure of radial solutions. Among other things, they proved the following result:

\medskip

\noindent
{\bf Theorem A} If $(u_1,u_2) $ is a radially symmetric non-topological solutions to the rank $2$ Chern-Simons system of Lie type with all vortices at the origin, then
\begin{equation}\label{infinity}
u_1(r)=-2\alpha_1\log r+O(1),\ u_2(r)=-2\alpha_2\log r+O(1)
\end{equation}
at infinity for some $\alpha_1,\alpha_2>1$. Futhermore,
\begin{equation*}
J(\alpha_1-1,\alpha_2-1)>J(N_1+1,N_2+1),
\end{equation*}
where $J(x,y)$ is the quadratic form associated to $K^{-1}$.

\medskip

For the existence of non-topological solutions of radially symmetric solutions, recently, in \cite{CKL}, Choe, Kim and the second author proved the following result:

\medskip

\noindent
{\bf Theorem B}
If $(\alpha_1,\alpha_2)$ defined in (\ref{infinity}) satisfies
\begin{eqnarray*}
&&-2N_1-N_2-3<\alpha_2-\alpha_1<2N_2+N_1+3,\\
&&2\alpha_1+\alpha_2>N_1+2N_2+6 \ \mbox{ and }\ \alpha_1+2\alpha_2>2N_1+N_2+6,
\end{eqnarray*}
then the ${\bf A}_2$ Chern-Simons system has a radially symmetric solution $(u_1,u_2) $ subject to the boundary condition (\ref{infinity}).

For general configuration vortices in $\R^2$, we first got the existence of non-topological solutions for the ${\bf A}_2$ and ${\bf B}_2$ Chern-Simons system by perturbation from the ${\bf A}_2$ and ${\bf B_2}$ Toda system with a singular source. In \cite{ALW}, we proved the following:

\medskip

\noindent
{\bf Theorem C}
Let $\{p_j\}_{j=1}^{N_1}$, $\{q_j\}_{j=1}^{N_2}\subset \R^2$. If either

(a) $N_2\displaystyle\sum_{j=1}^{N_1} p_j=N_1\displaystyle \sum_{j=1}^{N_2} q_j $ ;

\noindent or

(b)  $N_2\displaystyle\sum_{j=1}^{N_1}p_j \not =N_1\displaystyle \sum_{j=1}^{N_2} q_j$ and  $N_1,N_2>1$, $
 |N_1-N_2|\neq 1$,
\noindent
then there exists a non-topological solution $(u_1,u_2)$ of the ${\bf A}_2$ and ${\bf B}_2$ Chern-Simons system respectively.

\subsection{Main Results} In this paper, we continue our work on the rank $2$ Chern-Simons system. We consider the remaining ${\bf G}_2$ Chern-Simons system. We give an affirmative answer to the existence of non-topological solutions for the system with Cartan  matrix ${\bf G}_2$. Our main theorem can be stated as follows.

\begin{theorem}\label{thm101}
Let $\{p_j\}_{j=1}^{N_1}$, $\{q_j\}_{j=1}^{N_2}\subset \R^2$. If either

(a) $N_2\displaystyle\sum_{j=1}^{N_1} p_j=N_1\displaystyle \sum_{j=1}^{N_2} q_j $ ;

\noindent or

(b)  $N_2\displaystyle\sum_{j=1}^{N_1}p_j \not =N_1\displaystyle \sum_{j=1}^{N_2} q_j$ and  $N_1,N_2>1$, $
 |N_1-N_2|\neq 1$,
\noindent
then there exists a non-topological solution $(u_1,u_2)$ of problem (\ref{eqofg}).
\end{theorem}

\medskip

\begin{remark}
Note that if $N_1=0$ or $N_2=0$, by translation,  assumption $(a)$ in Theorem \ref{thm101} is always satisfied.
\end{remark}

\medskip
Non-topological solutions play very important role in the bubbling analysis of solutions to \eqref{equation1}. Therefore, our result is only the first step towards understanding the solution structure of non-topological solution of \eqref{equation1}. For further study on non-topological solutions for the Abelian case, we refer to \cite{cfl} and \cite{ckl}.

\medskip

We will prove Theorem \ref{thm101} in three cases which we describe below:
\begin{eqnarray}
\mbox{ {\bf Assumption (i)}:}&&\  \sum_{j=1}^{N_1} p_j= \sum_{j=1}^{N_2} q_j \ ,\ N_1=N_2;\\
\mbox{ {\bf Assumption (ii)}:}&& \ N_2\sum_{j=1}^{N_1} p_j= N_1\sum_{j=1}^{N_2} q_j \ , \ N_1\not =N_2\ , \ N_1, N_2\neq 1 \nonumber\\
&&\mbox{ or }N_2\sum_{j=1}^{N_1} p_j\neq N_1\sum_{j=1}^{N_2} q_j\ , \ |N_1-N_2|\neq 1 \ , \ N_1,N_2>1;\\
\mbox{{\bf Assumption (iii)}:}&& N_2\sum_{j=1}^{N_1} p_j= N_1\sum_{j=1}^{N_2} q_j,\ N_1\neq N_2, \ N_1=1 \mbox{ or } N_2=1.
\end{eqnarray}
If we can prove the existence of non-topological solutions under the above three assumptions separately, then it is easy to see that Theorem \ref{thm101} is proved. So in the following, we will prove the theorem under the three assumptions respectively.

\subsection{Sketch of the Proof }

In the following, we will outline the sketch of our proof. We follow exactly the same idea as in the proof for the ${\bf A}_2$ and ${\bf B}_2$ case.

\medskip

As in \cite{ALW}, we will view equation \eqref{eqofg} as a small perturbation of the ${\bf G}_2$ Toda system with a singular source.

After a suitable scaling transformation, the system (\ref{eqofg}) is transformed to

\begin{equation}\label{3m}
\left\{\begin{array}{c}
\Delta \tilde{U}_1+\Pi_{j=1}^{N_1}|\tilde{z}-\ve p_j|^2e^{2\tilde{U}_1-\tilde{U}_2}\\
=2\ve^2\Pi_{j=1}^{N_1}|\tilde{z}-\ve p_j|^4e^{4\tilde{U}_1-2\tilde{U}_2}
-\ve^2\Pi_{j=1}^{N_1}|\tilde{z}-\ve p_j|^2\Pi_{j=1}^{N_2}|\tilde{z}-\ve q_j|^2e^{\tilde{U}_2-\tilde{U}_1},\\
\mbox{$\ \ \ \ \ \ $} \\
\Delta \tilde{U}_2+\Pi_{j=1}^{N_2}|\tilde{z}-\ve q_j|^2e^{2\tilde{U}_2-3\tilde{U}_1}\\
=2\ve^2\Pi_{j=1}^{N_2}|\tilde{z}-\ve q_j|^4e^{4\tilde{U}_2-6\tilde{U}_1}
-3\ve^2\Pi_{j=1}^{N_2}|\tilde{z}-\ve q_j|^2\Pi_{j=1}^{N_1}|\tilde{z}-\ve p_j|^2e^{\tilde{U}_2-\tilde{U}_1}.
\end{array}
\right.
\end{equation}
When $\ve=0$, we obtain the following limiting system
\begin{equation}
\label{limitg}
\left\{\begin{array}{l}
\Delta \tilde{U}_1+|z|^{2N_1} e^{2\tilde{U}_1- \tilde{U}_2}=0  \ \ \mbox{in} \ \R^2 \\
\Delta \tilde{U}_2+ |z|^{2N_2} e^{2\tilde{U}_2- 3\tilde{U}_1}=0 \ \ \mbox{in} \ \R^2 \\
\int_{\R^2} |z|^{2N_1} e^{2\tilde{U}_1-\tilde{U}_2} <+\infty \ , \  \int_{\R^2} |z|^{2N_2} e^{2\tilde{U}_2-3\tilde{U}_1} <+\infty
\end{array}
\right.
\end{equation}
which is the  ${\bf G}_2$ Toda system with a single source at the origin.

\medskip

An immediate problem is the classification and non-degeneracy of the above system. In \cite{lwy}, Lin, Wei and Ye obtained the classification and non-degeneracy results of the $SU(N+1)$ Toda system with singular sources. In \cite{ALW1}, we use the results of \cite{lwy} to obtain a complete classification and non-degeneracy of the ${\bf G}_2$ Toda system (\ref{limitg}). In fact, the Toda system with ${\bf G}_2$ can be embedded into  the ${\bf A}_6$ Toda system under suitable group action. The solutions to (\ref{limitg}) depend on fourteen parameters
\begin{eqnarray*}
({\bf a},\lambda)&&=
(c_{43,1},c_{43,2},c_{52,1},c_{52,2},c_{53,1},c_{53,2},c_{54,1},c_{54,2},\\
&&c_{61,1},c_{61,2},c_{62,1},
c_{62,2}, \lambda_4,\lambda_5).
\end{eqnarray*}
The dimension of the linearized operator is {\bf fourteen}.

\medskip

The main difficulty of dealing with system is the large dimension of kernels. There are no explicit formula for the coefficients, except in  the case $ \sum_{j=1}^{N_1} p_j=\sum_{j=1}^{N_2} q_j \ , \ N_1=N_2$  which can be considered as  the reminiscent of the $SU(2)$ scalar equation. To get over this difficulty, we make use of the two scaling parameters for solutions of the Toda system and introduce two more free parameters . Instead of solving  the coefficient matrices for fixed scaling parameters, we only need to compute the two matrices in front of the two free scaling parameters we introduce.

\medskip

Now let us be more specific. The term of  order $O(\ve)$ will satisfy (\ref{Psi0b}) and (\ref{Psi1b}) in Section \ref{sec3.4}. The $O(\ve^2) $ term will satisfy (\ref{psib}). In this $O(\ve^2)$ term $\psi$, we introduce {\bf two free parameters $\xi_1,\xi_2$ }which play an important role in our proof. See Section \ref{sec3.6} and \ref{sec3.7}. At last the solution we find will be of this form
\begin{equation}\label{U}
\tilde{U}=\tilde{U}_{\mathbf{b}}+\ve\Psi+\ve^2 \psi+\ve^2v,
\end{equation}
where $\tilde{U}_{\mathbf{b}}=(\tilde{U}_{1,\mathbf{b}},\tilde{U}_{2,\mathbf{b}})$ is the solution of (\ref{limitg}), and  ${\bf b}$ denote the parameters $(\lambda, {\bf a})$ for simplicity of notations.
In order to solve in $v$, we need to solve a linearized problem:
\begin{equation}
\left\{\begin{array}{c}
\Delta \phi_1+|z|^{2N_1}e^{2\tilde{U}_{1,0}-\tilde{U}_{2,0}}(2\phi_1-\phi_2)=f_1\\
\Delta \phi_2+|z|^{2N_1}e^{2\tilde{U}_{2,0}-3\tilde{U}_{1,0}}(2\phi_2-3\phi_1)=f_2
\end{array}
\right.
\end{equation}
where $f_1$ and $f_2$ are explicitly given. Due to the existence of the kernel of the linearized equation, $(f_1,f_2)$ must satisfy some extra condition in order to have a solution. See Lemma \ref{lemma2b} for the necessary and sufficient condition. After that, we use the Liapunov-Schmidt reduction method to solve the nonlinear equation. It turns out that we can choose the perturbation $\mathbf{a}$ and $\lambda$ such that we can get the solution.

\medskip

Now we comment on the technical conditions. In the proof, we will choose $(\lambda_4,\lambda_5)$ first, depending on the assumptions. In general, the reduced problem for ${\bf a}$ has the following  form
\begin{equation}
\label{Q1010}
\frac{1}{\ve} {\bf B} {\bf a}+ \frac{1}{\ve} {\bf A} {\bf a} \cdot {\bf a}
+  {\bf Q} {\bf a} + O(|{\bf a}|^2) +{\bf a}_0  = O(\ve),
\end{equation}
where ${\bf A}, {\bf B}, {\bf Q}$ are matrices of size $12\times 12$, and $ {\bf a}_0 \in \R^{12}$.  Furthermore, the matrix ${\bf Q}$ can be decomposed into
\begin{equation}\label{qa}
{\bf Q}= \xi_1 {\bf Q}_1+ \xi_2 {\bf Q}_2 + {\mathcal T}
\end{equation}
where $ \xi_1$ and $ \xi_2$ are two free parameters. As we said before, we shall not attempt to compute the matrices ${\bf A}, {\bf B}$ and ${\mathcal T}$. Instead we focus on the two matrices ${\bf Q}_1$ and ${\bf Q}_2$. All we need to show is that at least of one of these two matrices is non-degenerate.

In case (a) of Theorem \ref{thm101}, i.e.  $ N_2\displaystyle\sum_{j=1}^{N_1} p_j=N_1\displaystyle\sum_{j=1}^{N_2} q_j$, by a shift of origin, we may assume that $ \displaystyle\sum_{j=1}^{N_1}p_j= \displaystyle\sum_{j=1}^{N_2} q_j =0$. In this case, the $\ve-$term $\ve \Psi$ vanishes and both ${\bf A}$ and ${\bf B}$ vanish. If $N_1\neq N_2, N_1, N_2\neq 1$ or $N_1=N_2$, then ${\bf a}_0$ vanishes and we obtain a reduced problems (in terms of ${\bf a}$) as follows:
\begin{equation}
\label{Q100}
 {\bf Q} {\bf a} + O(|{\bf a}|^2) = O(\ve).
\end{equation}
If $N_1\neq N_2, N_1=1$, or $N_2=1$, we use a different $O(\ve^2)$ approximation $\psi$ in (\ref{U}) and we obtain the reduced problem as follows:
\begin{equation}
{\bf Q}{\bf a}+O(|{\bf a}|^2)+{\bf a}_0=O(\ve).
\end{equation}
See Section \ref{sec3.7}. In both cases, we can show that the matrix ${\bf Q}$ is non-degenerate and (\ref{Q1010}) can be solved by contraction mapping.

The case (b) is considerably more difficult. Since $N_2\displaystyle\sum_{j=1}^{N_1}p_j \not = N_1 \displaystyle\sum_{j=1}^{N_2} q_j$, the $\ve-$ term exists and presents great difficulty in solving the reduced problem (\ref{Q1010}) in ${\bf a}$. To show that ${\bf B}=0$, we need $ |N_1-N_2| \not = 1$.  To show that the ${\bf a}_0$ term vanishes, we need $ N_1, N_2 >1$. In this case the reduced problem now takes the form
\begin{equation}
\label{Q101}
\frac{1}{\ve} {\bf A} {\bf a} \cdot {\bf a} +  {\bf Q} {\bf a} + O(|{\bf a}|^2) = O(\ve).
\end{equation}
where ${\bf Q}$ has the form of (\ref{qa}). By choosing large $\xi_1$ and $\xi_2=0$, we can solve (\ref{Q101}) such that $|{\bf a}|\leq O(\ve)$.

\medskip

In summary, the technical condition we have imposed is to make sure that ${\bf B}=0$ and that the quadratic term $\frac{1}{\ve} {\bf A} {\bf a} \cdot {\bf a} $ and the $O(1)$ term can not coexist.

\medskip

The organization of the paper is the following. In Section \ref{sec3}, we  present several important preliminaries of analysis. We first formulate our problem in terms of the functional equations (Section \ref{sec3.1}). Then we apply the classification and nondegeneracy result of ${\bf G}_2$ Toda system (Section \ref{sec3.2}).  In Section \ref{sec3.3}, we establish the invertibility properties of the linearized operator. Finally we obtain the next two orders  $O(\ve)$ and $O(\ve^2)$ in Section \ref{sec3.4}. In Section \ref{sec4}, we solve a projected nonlinear problem based on the preliminary results in Section 2. In Section \ref{sec3.5}, we prove our main theorem under the {\bf Assumption (i)}. In Section \ref{sec3.6}, we prove the theorem under the {\bf Assumption (ii)}. In Section \ref{sec3.7}, we prove the theorem under the {\bf Assumption (iii)}.

\section{Preliminaries}\label{sec3}
In this section, we consider the following ${\bf G}_2$ system in $\R^2$:
\begin{equation}\label{1b}
\left\{\begin{array}{c}
\Delta u_1+2e^{u_1}-e^{u_2}=4e^{2u_1}-2e^{2u_2}+e^{u_1+u_2}+4\pi\displaystyle\sum_{j=1}^{N_1}\delta_{p_j}\\
\Delta u_2+2e^{u_2}-3e^{u_1}=4e^{2u_2}-6e^{2u_1}-3e^{u_1+u_2}+4\pi\displaystyle\sum_{j=1}^{N_2}\delta_{q_j}.
\end{array}
\right.
\end{equation}

\subsection{Functional Formulation of the Problem }\label{sec3.1}

Defining
\begin{equation*}
u_1=\sum_{j=1}^{N_1}\ln|z-p_j|^2+\tilde{u}_1, \ \
u_2=\sum_{j=1}^{N_2}\ln|z-q_j|^2+\tilde{u}_2,
\end{equation*}
and $z=\frac{\tilde{z}}{\ve}$, and let $U_i$ and $\tilde{U}_i$ to be
$$
\tilde{u}_1(z)=U_1(\tilde{z})+(2N_1+2)\ln \ve, \ \ \tilde{u}_2(z)=U_2(\tilde{z})+(2N_2+2)\ln \ve,
$$
and
$$
\vect{\tilde{U}_1}{\tilde{U}_2}={\bf G}_2^{-1}\vect{U_1}{U_2},
$$
where ${\bf G}_2$ is the Cartan matrix $\left(\begin{array}{cc}
2&-1\\
-3&2
\end{array}
\right)$. Then  $ (\tilde{U}_1, \tilde{U}_2)$ will satisfy

\begin{equation}\label{3b}
\left\{\begin{array}{c}
\Delta \tilde{U}_1+\Pi_{j=1}^{N_1}|\tilde{z}-\ve p_j|^2e^{2\tilde{U}_1-\tilde{U}_2}\\
=2\ve^2\Pi_{j=1}^{N_1}|\tilde{z}-\ve p_j|^4e^{4\tilde{U}_1-2\tilde{U}_2}
-\ve^2\Pi_{j=1}^{N_1}|\tilde{z}-\ve p_j|^2\Pi_{j=1}^{N_2}|\tilde{z}-\ve q_j|^2e^{\tilde{U}_2-\tilde{U}_1},\\
\mbox{$\ \ \ \ \ \ $} \\
\Delta \tilde{U}_2+\Pi_{j=1}^{N_2}|\tilde{z}-\ve q_j|^2e^{2\tilde{U}_2-3\tilde{U}_1}\\
=2\ve^2\Pi_{j=1}^{N_2}|\tilde{z}-\ve q_j|^4e^{4\tilde{U}_2-6\tilde{U}_1}
-3\ve^2\Pi_{j=1}^{N_2}|\tilde{z}-\ve q_j|^2\Pi_{j=1}^{N_1}|\tilde{z}-\ve p_j|^2e^{\tilde{U}_2-\tilde{U}_1}.
\end{array}
\right.
\end{equation}

From now on, we shall work with (\ref{3b}). For simplicity of notations, we still denote the variable by $z$ instead of $\tilde{z}$.

\subsection{First Approximate Solution}\label{sec3.2}
When $\ve=0$, (\ref{3b}) becomes
\begin{equation}\label{firstapproximationb}
\left\{\begin{array}{c}
\Delta \tilde{U}_1+|z|^{2N_1}e^{2\tilde{U}_1-\tilde{U}_2}=0\\
\Delta \tilde{U}_2+|z|^{2N_2}e^{2\tilde{U}_2-3\tilde{U}_1}=0\\
\int_{\R^2}|z|^{2N_1}e^{2\tilde{U}_1-\tilde{U}_2}<\infty,\ \ \int_{\R^2}|z|^{2N_2}e^{2\tilde{U}_2-3\tilde{U}_1}<\infty
\end{array}
\right.
\end{equation}

For this system, we have gotten the classification and non-degeneracy result in \cite{ALW1}. From Theorem 2.1 in \cite{ALW1}, we see that all the solutions of (\ref{firstapproximationb}) depend on fourteen parameters $(c_{43},c_{52},c_{53},c_{54},c_{61},c_{62},\lambda_4,\lambda_5)\in \mathbb{C}^6\times(\mathbb{R}^{+})^2$, and all the solutions of (\ref{firstapproximationb}) are of the form
\begin{equation}
e^{-\tilde{U}_1}=2(\lambda_0+\sum_{i=1}^6\lambda_i|P_i(z)|^2),
\end{equation}
where
\begin{equation}
P_i(z)=z^{\mu_1+\cdots+\mu_i}+\sum_{j=0}^{i-1}c_{ij}z^{\mu_1+\cdots+\mu_j},
\end{equation}
\begin{eqnarray*}
\mu_1=\mu_3=\mu_4=\mu_6=N_1+1,\ \mu_2=\mu_5=N_2+1,
\end{eqnarray*}
and
\begin{eqnarray}
&&\lambda_0=\frac{1}{(2^{11-\frac{1}{2}}\mu_1^2\mu_2(\mu_1+\mu_2)^2(2\mu_1+\mu_2)^2(3\mu_1+\mu_2)
(3\mu_1+2\mu_2))^2\lambda_4\lambda_5},\\
&&\lambda_1=\frac{1}{(2^7 \mu_1\mu_2(\mu_1+\mu_2)^2(2\mu_1+\mu_2)(3\mu_1+2\mu_2))^2\lambda_5},\\
&&\lambda_2=\frac{1}{(2^7\mu_1^2\mu_2(\mu_1+\mu_2)(2\mu_1+\mu_2)(3\mu_1+\mu_2))^2\lambda_4},\\
&&\lambda_3=\frac{1}{(2^3\mu_1(\mu_1+\mu_2)(2\mu_1+\mu_2))^2},\\
&&\lambda_6=(2^{3+\frac{1}{2}}\mu_1\mu_2(\mu_1+\mu_2))^2\lambda_4\lambda_5.
\end{eqnarray}

See Theorem 2.1 in \cite{ALW1}.

We denote by
\begin{equation}
(\lambda,{\bf a})= (\lambda_4,\lambda_5,c_{43,1},c_{43,2},c_{52,1},c_{52,2},
 c_{53,1},c_{53,2},c_{54,1},c_{54,2},c_{61,1},c_{61,2},c_{62,1},c_{62,2}).
\end{equation}

When $\mathbf{a}=0$, the radially symmetric solution of (\ref{firstapproximationb}) can be expressed as follows:
\begin{eqnarray}
&&e^{-\tilde{U}_{1,0}}:=\rho_{1,G}^{-1}=2\rho_1^{-1},\\
&&e^{-\tilde{U}_{2,0}}:=\rho_{2,G}^{-1}=4\rho_2^{-1}
\end{eqnarray}
and
\begin{eqnarray}
\rho_{1}^{-1}&=&\lambda_0+\lambda_1r^{2\mu_1}+\lambda_2r^{2(\mu_1+\mu_2)}+\lambda_3r^{2(2\mu_1+\mu_2)}
+\lambda_4r^{2(3\mu_1+\mu_2)}+\lambda_5r^{2(3\mu_1+2\mu_2)}+\lambda_6r^{2(4\mu_1+2\mu_2)},\nonumber\\
\rho_{2}^{-1}&=&4\Big[\lambda_0 \mu_1^2 \lambda_1+r^{4 (\mu_1+\mu_2)} (4 r^{2 \mu_1} (\lambda_0 \lambda_6 (2 \mu_1+\mu_2)^2+\lambda_1 \lambda_5 (\mu_1+\mu_2)^2)+\lambda_0 \lambda_5 (3 \mu_1+2 \mu_2)^2\nonumber\\
&&+r^{4 \mu_1} (\lambda_1 \lambda_6 (3 \mu_1+2 \mu_2)^2+\mu_1^2 \lambda_3 \lambda_4)+\mu_1^2 \lambda_2 (\lambda_3+4 \lambda_4 r^{2 \mu_1}))\nonumber\\
&&+r^{2 \mu_2} (\lambda_0 (\lambda_2 (\mu_1+\mu_2)^2+\lambda_3 (2 \mu_1+\mu_2)^2 r^{2 \mu_1}+\lambda_4 (3 \mu_1+\mu_2)^2 r^{4 \mu_1})\nonumber\\
&&+\lambda_1 r^{2 \mu_1} (\mu_2^2 \lambda_2+\lambda_3 (\mu_1+\mu_2)^2 r^{2 \mu_1}+\lambda_4 (2 \mu_1+\mu_2)^2 r^{4 \mu_1}))\nonumber\\
&&+r^{6 (\mu_1+\mu_2)} (r^{2 \mu_1} (\lambda_2 \lambda_6 (3 \mu_1+\mu_2)^2+\lambda_3 \lambda_5 (\mu_1+\mu_2)^2)+\lambda_2 \lambda_5 (2 \mu_1+\mu_2)^2\nonumber\\
&&+r^{4 \mu_1} (\lambda_3 \lambda_6 (2 \mu_1+\mu_2)^2+\mu_2^2 \lambda_4 \lambda_5)+\lambda_4 \lambda_6 (\mu_1+\mu_2)^2 r^{6 \mu_1})+\mu_1^2 \lambda_5 \lambda_6 r^{12 \mu_1+8 \mu_2}\Big],
\end{eqnarray}
where $\mu_1=N_1+1, \ \ \mu_2=N_2+1$ and $\lambda_i$ are defined before.

Observe that the radial solution $ (\tilde{U}_{1,0}, \tilde{U}_{2,0})$ depends on two scaling parameters $(\lambda_4, \lambda_5)$. Later we shall choose $(\lambda_4,\lambda_5)$ in different ways.

\medskip

Next we have the following non-degeneracy result:
\begin{lemma}\label{lemma1b}
(Non-degeneracy) The previous solutions of (\ref{firstapproximationb}) are non-degenerate, i.e., the set of solutions corresponding to the linearized operator at $(\tilde{U}_{1,0},\tilde{U}_{2,0})$ is exactly fourteen dimensional. More precisely, if $\phi=\vect{\phi_1}{\phi_2}$ satisfies $|\phi (z)|\leq C(1+|z|)^{\alpha}$ for some $0\leq \alpha <1$, and
\begin{equation}\label{e210b}
\left\{\begin{array}{c}
\Delta \phi_1+|z|^{2N_1}e^{2\tilde{U}_{1,0}-\tilde{U}_{2,0}}(2\phi_1-\phi_2)=0\\
\Delta \phi_2+|z|^{2N_2}e^{2\tilde{U}_{2,0}-3\tilde{U}_{1,0}}(2\phi_2-3\phi_1)=0,
\end{array}
\right.
\end{equation}
then $\phi$ belongs to the following linear space $\mathcal{K}$: \ the span of
$$
\{
Z_{\lambda_4},Z_{\lambda_5},Z_{c_{43,1}},Z_{c_{43,2}},Z_{c_{52,1}},Z_{c_{52,2}},Z_{c_{53,1}},Z_{c_{53,2}},
Z_{c_{54,1}},Z_{c_{54,2}},Z_{c_{61,1}},Z_{c_{61,2}},Z_{c_{62,1}},Z_{c_{62,2}}
\},
$$
where
\begin{eqnarray*}
&&Z_{\lambda_4}=\vect{Z_{\lambda_4,1}}{Z_{\lambda_4,2}}=\vect{\partial_{\lambda_4}\tilde{U}_{1,0}}{
\partial_{\lambda_4}\tilde{U}_{2,0}},\ \ Z_{\lambda_5}=\vect{Z_{\lambda_5,1}}{Z_{\lambda_5,2}}=\vect{\partial_{\lambda_5}\tilde{U}_{1,0}}{
\partial_{\lambda_5}\tilde{U}_{2,0}},\\
&&Z_{c_{43,i}}=\vect{Z_{c_{43,i},1}}{Z_{c_{43,i},2}}=\vect{\partial_{c_{43,i}}\tilde{U}_{1,0}}{
\partial_{c_{43,i}}\tilde{U}_{2,0}},\ \ Z_{c_{52,i}}=\vect{Z_{c_{52,i},1}}{Z_{c_{52,i},2}}=\vect{\partial_{c_{52,i}}\tilde{U}_{1,0}}{
\partial_{c_{52,i}}\tilde{U}_{2,0}},\\
&&Z_{c_{53,i}}=\vect{Z_{c_{53,i},1}}{Z_{c_{53,i},2}}=\vect{\partial_{c_{53,i}}\tilde{U}_{1,0}}{
\partial_{c_{53,i}}\tilde{U}_{2,0}},\ \ Z_{c_{54,i}}=\vect{Z_{c_{54,i},1}}{Z_{c_{54,i},2}}=\vect{\partial_{c_{54,i}}\tilde{U}_{1,0}}{
\partial_{c_{54,i}}\tilde{U}_{2,0}},\\
&&Z_{c_{61,i}}=\vect{Z_{c_{61,i},1}}{Z_{c_{61,i},2}}=\vect{\partial_{c_{61,i}}\tilde{U}_{1,0}}{
\partial_{c_{61,i}}\tilde{U}_{2,0}},\ \ Z_{c_{62,i}}=\vect{Z_{c_{62,i},1}}{Z_{c_{62,i},2}}=\vect{\partial_{c_{62,i}}\tilde{U}_{1,0}}{
\partial_{c_{62,i}}\tilde{U}_{2,0}},
\end{eqnarray*}
for $i=1,2$.
\end{lemma}

For a proof, we refer Corollary 2.3 in \cite{ALW1}. For the reference of computation later, we want to write down the explicit expression of all kernel functions of the linearized equation. For the simplicity of notations, we use $\tilde{U}_{1,\lambda_4}, \cdots$ etc to denote $\partial_{\lambda_4}\tilde{U}_1, \cdots$ etc.
\begin{eqnarray*}
\tilde{U}_{1,\lambda_4}&=&\rho_{1}\Big(r^{2 (3 \mu_1 + \mu_2)} +2^7 \lambda_5 r^{4 (2 \mu_1 + \mu_2)} \mu_1^2 \mu_2^2 (\mu_1 + \mu_2)^2 -\frac{ r^{2 (\mu_1 + \mu_2)}}{2^{14} \lambda_4^2 \mu_1^4 \mu_2^2 (\mu_1 + \mu_2)^2 (2 \mu_1 +\mu_2)^2 (3 \mu_1 + \mu_2)^2} \\
&&- \frac{1}{2^{21} \lambda_4^2 \lambda_5 \mu_1^4 \mu_2^2 (\mu_1 + \mu_2)^4 (2 \mu_1 +\mu_2)^4 (3 \mu_1 + \mu_2)^2 (3 \mu_1 + 2 \mu_2)^2}\Big),\\
\tilde{U}_{2,\lambda_4}&=&\frac{4\rho_{2}}{2^{25} \mu_2^4 (\mu_1+\mu_2)^8}\Big[-\frac{1}{\lambda_4^2 \lambda_5^2 (2 \mu_1+\mu_2)^6 (3 \mu_1+\mu_2)^2 (3 \mu_1^2+2 \mu_1 \mu_2)^4}\\
&&+\frac{2^{21} \mu_2^2 \lambda_5 (\mu_1+\mu_2)^6 r^{6 (\mu_1+\mu_2)} }{\mu_1^4 \lambda_4^2 (3 \mu_1+\mu_2)^2}\\
&&\ \ \times \Big(2^{14} \mu_1^4 \mu_2^4 \lambda_4^2 (\mu_1+\mu_2)^2 (3 \mu_1+\mu_2)^2 r^{4 \mu_1} (256 \mu_1^2 \lambda_4 (\mu_1+\mu_2)^4 r^{2 \mu_1}+3)-1\Big)\\
&&+\frac{3\times 2^{14} \mu_2^2 (\mu_1+\mu_2)^4 r^{4 (\mu_1+\mu_2)} (2^{14} \mu_1^4 \mu_2^2 \lambda_4^2 (\mu_1+\mu_2)^2 (2 \mu_1+\mu_2)^2 (3 \mu_1+\mu_2)^2 r^{4 \mu_1}-1)}{\mu_1^4 \lambda_4^2 (2 \mu_1+\mu_2)^4 (3 \mu_1+\mu_2)^2}\\
&&+\frac{2 (\mu_1+\mu_2)^2 r^{2 \mu_2} (-\frac{(\mu_1+\mu_2)^2}{\lambda_4^3 (2 \mu_1+\mu_2)^6 (3 \mu_1+\mu_2)^4}-\frac{192 \mu_1^2 \mu_2^2 r^{2 \mu_1}}{\lambda_4^2 (2 \mu_1+\mu_2)^4 (3 \mu_1+\mu_2)^2}+2^{20} \mu_1^6 \mu_2^2 (\mu_1+\mu_2)^2 r^{6 \mu_1})}{\mu_1^8 \lambda_5 (3 \mu_1+2 \mu_2)^2}\\
&&+2^{42} \mu_1^4 \mu_2^6 \lambda_5^2 (\mu_1+\mu_2)^{10} r^{12 \mu_1+8 \mu_2}\Big],
\end{eqnarray*}
\begin{eqnarray*}
\tilde{U}_{1,\lambda_5}&=&\rho_{1}\Big[r^{6 \mu_1 + 4 \mu_2} +2^7 \lambda_4 r^{4 (2 \mu_1 + \mu_2)} \mu_1^2 \mu_2^2 (\mu_1 + \mu_2)^2- \frac{r^{2 \mu_1}}{2^{14} \lambda_5^2 \mu_1^2 \mu_2^2 (\mu_1 + \mu_2)^4 (2 \mu_1 + \mu_2)^2 (3 \mu_1 + 2 \mu_2)^2}\\
&&- \frac{1}{2^{21} \lambda_4 \lambda_5^2 \mu_1^4 \mu_2^2 (\mu_1 + \mu_2)^4 (2 \mu_1 +\mu_2)^4 (3 \mu_1 + \mu_2)^2 (3 \mu_1 + 2 \mu_2)^2}\Big],\\
\tilde{U}_{2,\lambda_5}
&=&4\rho_{2}\Big[3 \mu_2^2 \lambda_4 r^{10 \mu_1+6 \mu_2}+2^8 \mu_1^4 \mu_2^2 \lambda_4 \lambda_5 (\mu_1+\mu_2)^2 r^{12 \mu_1+8 \mu_2}+2^7 \mu_1^2 \mu_2^2 \lambda_4^2 (\mu_1+\mu_2)^4 r^{6 (2 \mu_1+\mu_2)}\\
&&+\frac{3 r^{8 \mu_1+6 \mu_2}}{2^7 \mu_1^2 (2 \mu_1+\mu_2)^2}-\frac{\lambda_4 r^{2 (3 \mu_1+\mu_2)}}{2^{14} \mu_1^2 \mu_2^2 \lambda_5^2 (\mu_1+\mu_2)^4 (3 \mu_1+2 \mu_2)^2}\\
&&-\frac{3 r^{2 (2 \mu_1+\mu_2)}}{2^{21} \mu_1^4 \mu_2^2 \lambda_5^2 (\mu_1+\mu_2)^4 (2 \mu_1+\mu_2)^4 (3 \mu_1+2 \mu_2)^2}\\
&&-\frac{r^{2 \mu_2}}{2^{35} \mu_1^8 \mu_2^4 \lambda_4^2 \lambda_5^2 (\mu_1+\mu_2)^4 (2 \mu_1+\mu_2)^6 (3 \mu_1+\mu_2)^4 (3 \mu_1+2 \mu_2)^2}\\
&&-\frac{1}{2^{34} \mu_1^6 \mu_2^4 \lambda_4 (\mu_1+\mu_2)^8 (3 \mu_1+\mu_2)^2}\Big(\frac{\mu_1^2}{\lambda_5^3 (2 \mu_1+\mu_2)^6 (3 \mu_1+2 \mu_2)^4}+\frac{3\ 2^6 \mu_2^2 (\mu_1+\mu_2)^2 r^{2 (\mu_1+\mu_2)}}{\lambda_5^2 (2 \mu_1+\mu_2)^4 (3 \mu_1+2 \mu_2)^2}\\
&&-2^{20} \mu_1^2 \mu_2^2 (\mu_1+\mu_2)^6 r^{6 (\mu_1+\mu_2)}\Big)\Big],
\end{eqnarray*}

\begin{eqnarray*}
\tilde{U}_{1,c_{43,1}}&=&\rho_{1}\Big(\lambda_4 r^{5 \mu_1 + 2 \mu_2} + \frac{\lambda_6(2 \mu_1 + \mu_2)r^{7 \mu_1 + 4 \mu_2}}{2 \mu_2} +\frac{\lambda_1 (\mu_1+\mu_2) r^{\mu_1}}{2 (3\mu_1+\mu_2)}\\
&&+ \frac{\lambda_3(\mu_1 + \mu_2) (2 \mu_1 + \mu_2)r^{3 \mu_1 + 2 \mu_2}}{2 \mu_2 (3 \mu_1 + \mu_2)}\Big)\cos\mu_1\theta,\\
\tilde{U}_{2,c_{43,1}}&=&\frac{4\rho_{2}r^{\mu_1+2 \mu_2}}{2 \mu_2 (3 \mu_1+\mu_2)}
\Big[(\mu_1+\mu_2)^2 (\lambda_0 \lambda_3 (2 \mu_1+\mu_2)^2+\lambda_1 \mu_2^2 \lambda_2)\\
&&+r^{4 \mu_1} \Big(2 r^{2 \mu_2} (\lambda_0 \lambda_6 (3 \mu_1+\mu_2) (3 \mu_1+2 \mu_2) (2 \mu_1+\mu_2)^2\\
&&+\mu_2 (\lambda_1 \lambda_5 (\mu_1+\mu_2)^2 (3 \mu_1+2 \mu_2)+2 \mu_1^2 \lambda_2 \lambda_4 (3 \mu_1+\mu_2)))\\
&&+3 \lambda_1 \mu_2 \lambda_4 (\mu_1+\mu_2) (2 \mu_1+\mu_2) (3 \mu_1+\mu_2)\Big)\\
&&+2 \mu_2 (2 \mu_1+\mu_2) r^{2 \mu_1} (\lambda_0 \lambda_4 (3 \mu_1+\mu_2)^2+\lambda_1 \lambda_3 (\mu_1+\mu_2)^2)\\
&&+(2 \mu_1+\mu_2) r^{2 (3 \mu_1+\mu_2)} \Big(2 (\mu_1+\mu_2) (\lambda_1 \lambda_6 (3 \mu_1+2 \mu_2)^2\\
&&+\mu_1^2 \lambda_3 \lambda_4)+(2 \mu_1+\mu_2) r^{2 \mu_2} (\lambda_2 \lambda_6 (3 \mu_1+\mu_2)^2+\lambda_3 \lambda_5 (\mu_1+\mu_2)^2)\Big)\\
&&+2 (\mu_1+\mu_2) (3 \mu_1+\mu_2) r^{4 (2 \mu_1+\mu_2)} (\lambda_3 \lambda_6 (2 \mu_1+\mu_2)^2+\mu_2^2 \lambda_4 \lambda_5)\\
&&+3 \mu_2 \lambda_4 \lambda_6 (\mu_1+\mu_2) (2 \mu_1+\mu_2) (3 \mu_1+\mu_2) r^{10 \mu_1+4 \mu_2}\Big]\cos\mu_1\theta,
\end{eqnarray*}

\begin{eqnarray*}
\tilde{U}_{1,c_{52,1}}
&=&\frac{\rho_{1}r^{2 \mu_1+\mu_2}}{(\mu_1+\mu_2) (3 \mu_1+2 \mu_2)}\Big[ (r^{2 \mu_1} (\mu_1 \lambda_4 (3 \mu_1+\mu_2)+(3 \mu_1+2 \mu_2) r^{2 \mu_2} (\lambda_5 (\mu_1+\mu_2)\\
&&-2 \lambda_6 (3 \mu_1+\mu_2) r^{2 \mu_1}))-\mu_1 \lambda_3 (\mu_1+\mu_2))\Big]\cos(2\mu_1+\mu_2)\theta,\\
\tilde{U}_{2,c_{52,1}}
&=&\frac{4\rho_{2}r^{2 \mu_1+\mu_2}}{(\mu_1+\mu_2) (3 \mu_1+2 \mu_2)}
\Big[r^{2 \mu_2}\Big (-2 (3 \mu_1+2 \mu_2) r^{2 \mu_1} (2 \lambda_0 \lambda_6 (2 \mu_1+\mu_2)^2 (3 \mu_1+\mu_2)\\
&&-\lambda_1 \mu_2 \lambda_5 (\mu_1+\mu_2)^2)+\lambda_0 \lambda_5 (\mu_1+\mu_2)^2 (3 \mu_1+2 \mu_2)^2\\
&&-2 (\mu_1+\mu_2) (3 \mu_1+\mu_2) r^{4 \mu_1} (\lambda_1 \lambda_6 (3 \mu_1+2 \mu_2)^2+\mu_1^2 \lambda_3 \lambda_4)\\
&&+\mu_1^2 \lambda_2 (\lambda_3 (\mu_1+\mu_2)^2-2 \mu_2 \lambda_4 (3 \mu_1+\mu_2) r^{2 \mu_1})\Big)
+\mu_1^2 (\lambda_0 \lambda_4 (3 \mu_1+\mu_2)^2+\lambda_1 \lambda_3 (\mu_1+\mu_2)^2)\\
&&+\mu_1 r^{4 (\mu_1+\mu_2)} \Big(-2 (3 \mu_1+2 \mu_2) (\lambda_2 \lambda_6 (3 \mu_1+\mu_2)^2+\lambda_3 \lambda_5 (\mu_1+\mu_2)^2)\\
&&-2 (\mu_1+\mu_2) r^{2 \mu_1} (\lambda_3 \lambda_6 (2 \mu_1+\mu_2)^2+\mu_2^2 \lambda_4 \lambda_5)+3 \lambda_4 \lambda_6 (\mu_1+\mu_2) (3 \mu_1+\mu_2) (3 \mu_1+2 \mu_2) r^{4 \mu_1}\Big)\\
&&+3 \mu_1 \lambda_5 \lambda_6 (\mu_1+\mu_2) (3 \mu_1+\mu_2) (3 \mu_1+2 \mu_2) r^{8 \mu_1+6 \mu_2}\Big]\cos(2\mu_1+\mu_2)\theta,
\end{eqnarray*}

\begin{eqnarray*}
\tilde{U}_{1,c_{53,1}}
&=&\rho_{1}\Big[\frac{\mu_1 \lambda_2 r^{\mu_1+\mu_2}}{2 (3 \mu_1+2 \mu_2)}-\frac{\mu_1 \lambda_3 (2 \mu_1+\mu_2) r^{3 \mu_1+\mu_2}}{2 \mu_2 (3 \mu_1+2 \mu_2)}+\lambda_5 r^{5 \mu_1+3 \mu_2}\\
&&-\frac{\lambda_6 (2 \mu_1+\mu_2) r^{7 \mu_1+3 \mu_2}}{2 \mu_2}\Big]\cos(\mu_1+\mu_2)\theta,\\
\tilde{U}_{2,c_{53,1}}
&=&\frac{4\rho_{2}}{2 \mu_2 (3 \mu_1+2 \mu_2)} \Big[r^{3 (\mu_1+\mu_2)} \Big(2 r^{2 \mu_1} (\mu_2 (2 \lambda_1 \lambda_5 (\mu_1+\mu_2)^2 (3 \mu_1+2 \mu_2)+\mu_1^2 \lambda_2 \lambda_4 (3 \mu_1+\mu_2))\\
&&-\lambda_0 \lambda_6 (2 \mu_1+\mu_2)^2 (3 \mu_1+\mu_2) (3 \mu_1+2 \mu_2))+2 \mu_2 (2 \mu_1+\mu_2) (\lambda_0 \lambda_5 (3 \mu_1+2 \mu_2)^2+\mu_1^2 \lambda_2 \lambda_3)\\
&&-(2 \mu_1+\mu_2)^2 r^{4 \mu_1} (\lambda_1 \lambda_6 (3 \mu_1+2 \mu_2)^2+\mu_1^2 \lambda_3 \lambda_4)\Big)-\mu_1^2 r^{\mu_1+\mu_2} (\lambda_0 \lambda_3 (2 \mu_1+\mu_2)^2+\lambda_1 \mu_2^2 \lambda_2)\\
&&-\mu_1 r^{5 (\mu_1+\mu_2)} (2 (2 \mu_1+\mu_2) r^{2 \mu_1} (\lambda_2 \lambda_6 (3 \mu_1+\mu_2)^2+\lambda_3 \lambda_5 (\mu_1+\mu_2)^2)\\
&&-3 \mu_2 \lambda_2 \lambda_5 (2 \mu_1+\mu_2) (3 \mu_1+2 \mu_2)+2 (3 \mu_1+2 \mu_2) r^{4 \mu_1} (\lambda_3 \lambda_6 (2 \mu_1+\mu_2)^2+\mu_2^2 \lambda_4 \lambda_5))\\
&&+3 \mu_1 \mu_2 \lambda_5 \lambda_6 (2 \mu_1+\mu_2) (3 \mu_1+2 \mu_2) r^{11 \mu_1+7 \mu_2}\Big]\cos(\mu_1+\mu_2)\theta,
\end{eqnarray*}

\begin{eqnarray*}
\tilde{U}_{1,c_{54,1}}
&=&\rho_{1}\Big[\lambda_5 r^{6 \mu_1 + 3 \mu_2} + \frac{\lambda_2 r^{2 \mu_1 + \mu_2} \mu_1 (3 \mu_1 + \mu_2)}{(\mu_1 + \mu_2)(3 \mu_1 + 2 \mu_2)}\Big]\cos\mu_2\theta,\\
\tilde{U}_{2,c_{54,1}}
&=&\frac{4\rho_{2}r^{\mu_2}}{(\mu_1+\mu_2) (3 \mu_1+2 \mu_2)}\Big[\lambda_0 (\mu_1+\mu_2) (3 \mu_1+\mu_2) (\mu_1^2 \lambda_2+\lambda_5 (3 \mu_1+2 \mu_2)^2 r^{2 (2 \mu_1+\mu_2)})\\
&&+r^{2 (2 \mu_1+\mu_2)} \Big(\lambda_5 (\mu_1+\mu_2)^2 (3 \mu_1+2 \mu_2) r^{2 \mu_1} (2 \lambda_1 (2 \mu_1+\mu_2)+\mu_1 r^{2 (\mu_1+\mu_2)} (\lambda_3+\lambda_6 r^{2 (2 \mu_1+\mu_2)}))\\
&&+\mu_1 \lambda_2 (\mu_1 (3 \mu_1+\mu_2) (\lambda_3 (\mu_1+\mu_2)+2 \lambda_4 (2 \mu_1+\mu_2) r^{2 \mu_1})\\
&&+r^{2 (\mu_1+\mu_2)} (6 \lambda_5 (\mu_1+\mu_2) (2 \mu_1+\mu_2)^2+\lambda_6 (3 \mu_1+\mu_2)^2 (3 \mu_1+2 \mu_2) r^{2 \mu_1}))\Big)\Big]\cos\mu_2\theta,
\end{eqnarray*}

\begin{eqnarray*}
\tilde{U}_{1,c_{61,1}}
&=&\rho_{1}\Big[\lambda_6 r^{5 \mu_1 + 2 \mu_2} + \frac{\lambda_5 r^{3 \mu_1 +2 \mu_2} \mu_2 (\mu_1 + \mu_2)}{(2 \mu_1 + \mu_2) (3 \mu_1 + \mu_2)}\Big]\cos(3\mu_1+2\mu_2)\theta,\\
\tilde{U}_{2,c_{61,1}}&=&\frac{4\rho_{2}r^{3 \mu_1+2 \mu_2}}{(2 \mu_1+\mu_2) (3 \mu_1+\mu_2)}\Big[2 \mu_1 (\lambda_0 \lambda_6 (2 \mu_1+\mu_2)^2 (3 \mu_1+\mu_2)-\lambda_1 \mu_2 \lambda_5 (\mu_1+\mu_2)^2)\\
&&-r^{2 \mu_2} \Big(\mu_2 \lambda_2 (2 \mu_1+\mu_2) (\lambda_5 (\mu_1+\mu_2)^2+\lambda_6 (3 \mu_1+\mu_2)^2 r^{2 \mu_1})\\
&&+\lambda_3 (\mu_1+\mu_2) (2 \mu_1+\mu_2) r^{2 \mu_1} (\mu_2 \lambda_5 (\mu_1+\mu_2)+\lambda_6 (2 \mu_1+\mu_2) (3 \mu_1+\mu_2) r^{2 \mu_1})\\
&&+\lambda_4 (\mu_1+\mu_2) (3 \mu_1+\mu_2) r^{4 \mu_1} (\mu_2^2 \lambda_5+\lambda_6 (2 \mu_1+\mu_2)^2 r^{2 \mu_1})\Big)\\
&&-6 \mu_1^2 \lambda_5 \lambda_6 (\mu_1+\mu_2) (2 \mu_1+\mu_2) r^{6 \mu_1+4 \mu_2}\Big]\cos(3\mu_1+2\mu_2)\theta,
\end{eqnarray*}

\begin{eqnarray*}
\tilde{U}_{1,c_{62,1}}
&=&\rho_{1}\Big[\lambda_6 r^{5 \mu_1 + 3 \mu_2} - \frac{\lambda_4 r^{3 \mu_1 + \mu_2} \mu_1 \mu_2}{(2 \mu_1+ \mu_2)(3 \mu_1 +2 \mu_2)}\Big]\cos(3\mu_1+\mu_2)\theta,\\
\tilde{U}_{2,c_{62,1}}
&=&\frac{4\rho_{2}r^{3 \mu_1+\mu_2}}{(2 \mu_1+\mu_2) (3 \mu_1+2 \mu_2)} \Big[r^{2 \mu_2} (2 (\mu_1+\mu_2) (\lambda_0 \lambda_6 (2 \mu_1+\mu_2)^2 (3 \mu_1+2 \mu_2)+\mu_1^2 \mu_2 \lambda_2 \lambda_4)\\
&&+\mu_2 (2 \mu_1+\mu_2) r^{2 \mu_1} (\lambda_1 \lambda_6 (3 \mu_1+2 \mu_2)^2+\mu_1^2 \lambda_3 \lambda_4))\\
&&+\mu_1^2 \lambda_1 \mu_2 \lambda_4 (2 \mu_1+\mu_2)-\mu_1 r^{4 (\mu_1+\mu_2)} ((3 \mu_1+2 \mu_2) (\lambda_3 \lambda_6 (2 \mu_1+\mu_2)^2+\mu_2^2 \lambda_4 \lambda_5)\\
&&+6 \lambda_4 \lambda_6 (\mu_1+\mu_2)^2 (2 \mu_1+\mu_2) r^{2 \mu_1})-\mu_1 \lambda_5 \lambda_6 (2 \mu_1+\mu_2)^2 (3 \mu_1+2 \mu_2) r^{6 (\mu_1+\mu_2)}\Big]\cos(3\mu_1+\mu_2)\theta,
\end{eqnarray*}

and by replacing the $\cos $ terms by $\sin$ terms, we have the $\tilde{U}_{c_{ij,2}}$.

\medskip

For simplicity of notations, we also denote by $(Z_1,Z_2,\cdots,Z_{14})$ the kernels $(Z_{\lambda_0},Z_{\lambda_1},\cdots,Z_{c_{62,2}})$.
Because $\{Z_i\}$ are linearly independent, we have

\begin{equation}\label{det1b}
\mbox{det}[(\int_{\R^2}\Delta Z_i\cdot Z_j)_{i,j=1,\cdots,14}]\neq 0.
\end{equation}

We have the following corollary:
\begin{corollary}\label{cor1}
If $\phi=\vect{\phi_1}{\phi_2}$ satisfies $|\phi (z)|\leq C(1+|z|)^{\alpha}$ for some $0\leq \alpha <1$, and
\begin{equation}\label{eq1}
\left\{\begin{array}{c}
\Delta \phi_1+2|z|^{2N_1}e^{2\tilde{U}_{1,0}-\tilde{U}_{2,0}}\phi_1-3|z|^{2N_2}e^{2\tilde{U}_{2,0}
-3\tilde{U}_{1,0}}\phi_2=0\\
\Delta \phi_2+2|z|^{2N_2}e^{2\tilde{U}_{2,0}-3\tilde{U}_{1,0}}\phi_2
-|z|^{2N_1}e^{2\tilde{U}_{1,0}-\tilde{U}_{2,0}}\phi_1=0,
\end{array}
\right.
\end{equation}
then $\phi$ belongs to the following linear space $\mathcal{K}^*$: \ the span of
$$
\{
Z_{\lambda_4}^*, Z_{\lambda_5}^*, Z_{c_{43,1}}^*, Z_{c_{43,2}}^*, Z_{c_{52,1}}^*, Z_{c_{52,2}}^*, Z_{c_{53,1}}^*, Z_{c_{53,2}}^*,Z_{c_{54,1}}^*,Z_{c_{54,2}}^*,Z_{c_{61,1}}^*,Z_{c_{61,2}}^*,Z_{c_{62,1}}^*,Z_{c_{62,2}}^*
\},
$$
where
\begin{equation}\label{222b}
Z_i^*=\vect{Z_{i,1}^*}{Z_{i,2}^*}=\vect{2Z_{i,1}-Z_{i,2}}{\frac{2}{3}Z_{i,2}-Z_{i,1}}.
\end{equation}
\end{corollary}

\medskip

We have
\begin{equation}\label{det2b}
\mbox{det}[(\int_{\R^2}Z_i^*\cdot Z_j^*)_{i,j=3,\cdots,14}]\neq 0.
\end{equation}

\medskip

We will choose the first approximate solution to be $\vect{\tilde{U}_{1,(\lambda,c_{43},c_{52},c_{53},c_{54},c_{61},c_{62})}}
{\tilde{U}_{2,(\lambda,c_{43},c_{52},c_{53},c_{54},c_{61},c_{62})}}$, where the parameters $\lambda,c_{43}, c_{52},c_{53},c_{54},c_{61},c_{62}$ satisfy
\begin{equation}\label{parab}
|\mathbf{a}|:=|c_{43}|+|c_{52}|+|c_{53}|+|c_{54}|+|c_{61}|+|c_{62}|\leq C_0\ve,\ \  |\lambda |=O(1)
\end{equation}
for some fixed constant $C_0>0$.

\medskip

For the simplicity of notations, we also denote ${\bf b}=(\lambda,{\bf a})$, and $\tilde{U}_{i,{\bf b}}=\tilde{U}_{i,(\lambda,c_{43},c_{52},c_{53},c_{54},c_{61},c_{62})}$. We want to look for solutions of the form
\begin{equation}
\tilde{U}_1=\tilde{U}_{1,{\bf b}}+\ve\Psi_{1}+\ve^2\phi_1,\ \
\tilde{U}_2=\tilde{U}_{2,{\bf b}}+\ve\Psi_{2}+\ve^2\phi_2,
\end{equation}
where ${\bf b}$ is fixed.

\medskip

To obtain the next order term, we need to  study the linearized operator around the solution  $\vect{\tilde{U}_{1,0}}{\tilde{U}_{2,0}}$.

\subsection{Invertibility of the Linearized Operator}\label{sec3.3}
Now we consider the invertibility of the linearized operator in some suitable Sobolev spaces. To this end, we use the technical framework introduced by Chae-Imanuvilov \cite{ci} and which has been used in \cite{ALW}. Let $\alpha \in (0,1)$ and
\begin{equation}
\label{X}
X_\alpha= \{ u \in L^2_{loc} (\R^2), \int_{\R^2} (1+|x|^{2+\alpha}) u^2 dx <+\infty \},
\end{equation}
\begin{equation}
\label{Y}
Y_\alpha=\{ u \in W^{2,2}_{loc} (\R^2), \int_{\R^2} (1+|x|^{2+\alpha}) | \Delta u|^2 + \frac{u^2}{1+|x|^{2+\alpha}} <+\infty \}.
\end{equation}

On $X_\alpha$ and $Y_\alpha$, we equip with two norms respectively:
\begin{equation}
\label{norm12}
\| f \|_{**}= \sup_{y \in \R^2} (1+|y|)^{2+\alpha}|f(y)|, \ \ \| h \|_{*}=\sup_{y \in \R^2} (\log (2+|y|))^{-1} |h(y)|.
\end{equation}

Clearly, the linearized operator in (\ref{e210b}) is bounded from $Y_\alpha$ to $X_\alpha$.

For $f=\vect{f_1}{f_2}$, $g=\vect{g_1}{g_2}$, we denote by
$\langle f,g\rangle=\int_{\R^2}f\cdot g dx$.

Using the non-degeneracy result we get, i.e. Lemma \ref{lemma1b} and Corollary \ref{cor1} in Section \ref{sec3.2} and noting that (\ref{eq1}) is the adjoint operator of (\ref{e210b}), we have the following:
\begin{lemma}
\label{lemma2b}
Assume that $h=\vect{h_1}{h_2} \in X_\alpha$ be such that
\begin{equation}
\langle Z_{i}^*,h\rangle=0, \mbox{ for } i=3,\cdots,14.
\end{equation}
Then  one can find a unique solution $ \phi=\vect{\phi_1}{\phi_2}=T^{-1}(h) \in Y_\alpha$ satisfying
\begin{equation}\label{e305b}
\left\{\begin{array}{c}
\Delta \phi_1+|z|^{2N_1}e^{2\tilde{U}_{1,0}-\tilde{U}_{2,0}}(2\phi_1-\phi_2)=h_1\\
\Delta \phi_2+|z|^{2N_1}e^{2\tilde{U}_{2,0}-3\tilde{U}_{1,0}}(2\phi_2-3\phi_1)=h_2
\end{array},
\right.
\| \phi \|_{*} \leq C \| h \|_{**},
\end{equation}
such that $<\Delta Z_i, \phi>=0$ for $i=1,\cdots, 14$.
Moreover, the map $ h \stackrel{T} {\longrightarrow} \phi$ can be made continuous and smooth .
\end{lemma}

We note that the uniqueness in Lemma \ref{lemma2b} is due to (\ref{det1b}). In the next subsections, we will use Lemma \ref{lemma2b} to obtain our approximate solution up to $O(\ve^2)$. However our approximation solution would be chosen according to our assumption of the vortex configuration.

\subsection{Improvements of the Approximate Solution}\label{sec3.4}
Similar to the ${\bf A}_2$ and ${\bf B}_2$ case, we need to find the $O(\ve) $ and $O(\ve^2)$ improvement of our approximate solutions. So we need to find solutions of the following equations.

Denote by
\begin{equation*}
f(\ve,z)=\Pi_1^{N_1}|z-\ve p_i|^2, \ \ g(\ve,z)=\Pi_1^{N_2}|z-\ve q_i|^2.
\end{equation*}
Then by Taylor's expansion, we have
$$ f(\ve, z)= f(0, z) + \ve f_{\ve} (0, z) + \frac{\ve^2}{2} f_{\ve \ve} (0, z) + O(\ve^3)$$
where
\begin{equation}\label{eq2}
 f(0, z)= |z|^{2N_1}, \ \  f_\ve (0, z)= - 2 |z|^{2N_1-2} <\sum_{j=1}^{N_1} p_j, z>.
\end{equation}
Similarly we can get the expansions for $ g(\ve, z)$.

Let $\vect{\Psi_{0,1}}{\Psi_{0,2}}$ be the solution of
\begin{equation}\label{Psi0b}
\left\{\begin{array}{c}
\Delta \Psi_{0,1}+|z|^{2N_1}e^{2\tilde{U}_{1,0}-\tilde{U}_{2,0}}(2\Psi_{0,1}-\Psi_{0,2})\\
\, \, \, \, =-f_{\ve}(0,z)e^{2\tilde{U}_{1,0}-\tilde{U}_{2,0}},\\
\Delta \Psi_{0,2}+|z|^{2N_2}e^{2\tilde{U}_{2,0}-3\tilde{U}_{1,0}}(2\Psi_{0,2}-3\Psi_{0,1})\\
\, \, \, \, =-g_{\ve}(0,z)e^{2\tilde{U}_{2,0}-3\tilde{U}_{1,0}}.
\end{array}
\right.
\end{equation}
Let $\vect{\Psi_{i,1}}{\Psi_{i,2}}$ be the solution of
\begin{equation}\label{Psi1b}
\left\{\begin{array}{c}
\Delta \Psi_{i,1}+|z|^{2N_1}e^{2\tilde{U}_{1,0}-\tilde{U}_{2,0}}(2\Psi_{i,1}-\Psi_{i,2})
=-|z|^{2N_1}e^{2\tilde{U}_{1,0}-\tilde{U}_{2,0}}(2\Psi_{0,1}-\Psi_{0,2})(2Z_{i,1}-Z_{i,2})\\
\ \ \ \ \ \ \ \ \ \ \ \ \ \ \ \ \ \ \ \ \ \ \ \ \ \ \ \ \ \ \ \ \ \ \ \ -f_\ve(0,z)e^{2\tilde{U}_{1,0}-\tilde{U}_{2,0}}(2Z_{i,1}-Z_{i,2})\\
\\
\Delta \Psi_{i,2}+|z|^{2N_2}e^{2\tilde{U}_{2,0}-3\tilde{U}_{1,0}}(2\Psi_{i,2}-3\Psi_{i,1})
=-|z|^{2N_2}e^{2\tilde{U}_{2,0}-3\tilde{U}_{1,0}}(2\Psi_{0,2}-3\Psi_{0,1})(2Z_{i,2}-3Z_{i,1})\\
\ \ \ \ \ \ \ \ \ \ \ \ \ \ \ \ \ \ \ \ \ \ \ \ \ \ \ \ \ \ \ \ \ \ \ \ \ \ -g_\ve(0,z)e^{2U_{2,0}-3U_{1,0}}(2Z_{i,2}-3Z_{i,1})
\end{array}
\right.
\end{equation}
for $i=3,\cdots,14$.

Let $\vect{\psi_1}{\psi_2}$ be the solution of

\begin{equation}\label{psib}
\left\{\begin{array}{c}
\Delta \psi_1+|z|^{2N_1}e^{2\tilde{U}_{1,0}-\tilde{U}_{2,0}}(2\psi_1-\psi_2)=
2|z|^{4N_1}e^{4\tilde{U}_{1,0}-2\tilde{U}_{2,0}}-|z|^{2(N_1+N_2)}e^{\tilde{U}_{2,0}-\tilde{U}_{1,0}}\\
-\frac{1}{2}|z|^{2N_1}e^{2U_{1,0}-U_{2,0}}(2\Psi_{0,1}-\Psi_{0,2})^2
-f_\ve(0,z)e^{2\tilde{U}_{1,0}-\tilde{U}_{2,0}}(2\Psi_{0,1}-\Psi_{0,2})
-\frac{f_{\ve\ve}(0,z)}{2}e^{2\tilde{U}_{1,0}-\tilde{U}_{2,0}}\\
\\
\Delta \psi_2+|z|^{2N_2}e^{2\tilde{U}_{2,0}-3\tilde{U}_{1,0}}(2\psi_2-3\psi_1)=
2|z|^{4N_2}e^{4\tilde{U}_{2,0}-6\tilde{U}_{1,0}}-3|z|^{2(N_1+N_2)}e^{\tilde{U}_{2,0}-\tilde{U}_{1,0}}\\
-\frac{1}{2}|z|^{2N_2}e^{2U_{2,0}-3U_{1,0}}(2\Psi_{0,2}-3\Psi_{0,1})^2
-g_\ve(0,z)e^{2\tilde{U}_{2,0}-3\tilde{U}_{1,0}}(2\Psi_{0,2}-3\Psi_{0,1})
-\frac{g_{\ve\ve}(0,z)}{2}e^{2\tilde{U}_{2,0}-3\tilde{U}_{1,0}}.
\end{array}
\right.
\end{equation}

\medskip

Obviously, if $\sum_{j=1}^{N_1}p_j=\sum_{j=1}^{N_2}q_j=0$, then $f_\ve(0,z)=g_\ve(0,z)=0$ by (\ref{eq2}), in this case, $\Psi_0=\Psi_i=0$ for $i=3,\cdots, 14$. Note that if $N_2\sum_{j=1}^{N_1}p_j=N_1\sum_{j=1}^{N_2}q_j$, we can always shift the origin such that $\sum_{j=1}^{N_1}p_j=\sum_{j=1}^{N_2}q_j=0$. Hence $f_\ve(0,z)\neq 0$ occurs only when assumption $(b)$ of Theorem \ref{thm101} holds. In this case, using Lemma \ref{lemma2b}, we know that there exists a unique solution $\vect{\psi_{0,1}}{\psi_{0,2}}\in Y_\alpha$ of (\ref{Psi0b}) such that $\langle \Psi_0,\Delta Z_j\rangle=0$ for $j=1,\cdots,14$.
If $N_2\sum_{j=1}^{N_1}p_j\neq N_1\sum_{j=1}^{N_2}q_j$, then $| N_1-N_2 | \neq 1$ implies that the right hand side of equation (\ref{Psi1b}) is orthogonal to $Z_i^*$ for $i=3,\cdots,14$.
By Lemma  \ref{lemma2b}, there exists a unique solution  $\vect{\Psi_{i,1}}{\Psi_{i,2}}\in Y_\alpha$ of (\ref{Psi1b}) such that $\langle \Psi_i,\Delta Z_j\rangle=0$ for $j=1,\cdots,14$. And if $N_1, N_2 >1$, the right hand side of (\ref{psib}) is  orthogonal to $Z_i^*$ for $i=3,\cdots,14$. By Lemma \ref{lemma2b}, we can find a unique solution $\psi_0=\vect{\psi_{0,1}}{\psi_{0,2}}\in Y_\alpha$ such that $\langle \psi_0, \Delta Z_i\rangle=0$ for $i=1,\cdots,14$.

\medskip

Similar to the ${\bf A}_2$ and ${\bf B}_2$ case, the solution we will use later  is $\psi=\psi_0+\xi_1Z_{\lambda_4}+\xi_2Z_{\lambda_5}$ where $\xi_1,\xi_2$ are two constants independent of ${\bf a}$ and will be determined later.

Finally, the approximate solution with all the terms of  $O(\ve)$ and $O(\ve^2)$ is
\begin{equation}\label{eq3}
V=\vect{V_1}{V_2}=\vect{\tilde{U}_{1,\bf{b}}+\ve(\Psi_{0,1}+\sum_{i=3}^{14}\Psi_{i,1}\mathbf{a}_i)+\ve^2\psi_1}
{\tilde{U}_{2,\bf{b}}+\ve(\Psi_{0,2}+\sum_{i=3}^{14}\Psi_{i,2}\mathbf{a}_i)+\ve^2\psi_2},
\end{equation}
we use the notation
\begin{eqnarray*}
\mathbf{b}&=&(\lambda_4,\lambda_5, \mathbf{a}),\\
\mathbf{a}&=&({\bf a}_3,{\bf a}_4,\cdots,{\bf a}_{14})\\
&=&(c_{43,1},c_{43,2},c_{52,1},c_{52,2},c_{53,1},c_{53,2},c_{54,1}, c_{54,2}, c_{61,1}, c_{61,2}, c_{62,1}, c_{62,2}).
\end{eqnarray*}

Parameters $\lambda_4,\lambda_5$ will be chosen later according to different assumptions in Theorem \ref{thm101}. After $\lambda_4,\lambda_5$ are fixed, ${\bf a}$ will be chosen in order to find a solution of (\ref{3b}) in the form of (\ref{eq3})
Then $\vect{V_1+\ve^2v_1}{V_2+\ve^2v_2}$ is a solution of (\ref{3b}) if $\vect{v_1}{v_2}$ satisfies
\begin{equation}\label{vb}
\left\{\begin{array}{c}
\Delta v_1+|z|^{2N_1}e^{2\tilde{U}_{1,0}-\tilde{U}_{2,0}}(2v_1-v_2)=G_1\\
\Delta v_2+|z|^{2N_2}e^{2\tilde{U}_{2,0}-3\tilde{U}_{1,0}}(2v_2-3v_1)=G_2,
\end{array}
\right.
\end{equation}
where
\begin{eqnarray}\label{eq4}
G_1&&= E_1+N_{11}(v)+N_{12}(v) ,        \\
G_2&&=E_2+N_{21}(v)+N_{22}(v),\label{eq5}\\
N_{11}(v)&&=2\Pi |z-\ve p_j|^4(e^{4\tilde{U}_1-2\tilde{U}_2}-e^{4V_1-2V_2})\nonumber\\
&&-\Pi |z-\ve p_j|^2\Pi |z-\ve q_j|^2(e^{\tilde{U}_2-\tilde{U}_1}-e^{V_2-V_1}),\nonumber\\
N_{12}(v)&&=\frac{-f(\ve,z)e^{2\tilde{U}_1-\tilde{U}_2}+f(\ve,z)e^{2V_1-V_2}
}{\ve^2}\nonumber\\
&&+f(0,z)e^{2\tilde{U}_{1,0}-\tilde{U}_{2,0}}(2v_1-v_2),\nonumber\\
N_{21}(v)&&=2\Pi |z-\ve q_j|^4(e^{4\tilde{U}_2-6\tilde{U}_1}-e^{4V_2-6V_1})\nonumber\\
&&-3\Pi |z-\ve p_j|^2\Pi |z-\ve q_j|^2(e^{\tilde{U}_2-\tilde{U}_1}-e^{V_2-V_1}),\nonumber\\
N_{22}(v)&&=\frac{-g(\ve,z)e^{2\tilde{U}_2
-3\tilde{U}_1}+g(\ve,z)e^{2V_2-3V_1}}{\ve^2}\nonumber\\
&&+g(0,z)e^{2\tilde{U}_{2,0}-3\tilde{U}_{1,0}}(2v_2-3v_1),\nonumber
\end{eqnarray}
and $E_i$ are the errors:
\begin{eqnarray*}
& \ \ & E_1 = 2\Pi |z-\ve p_j|^4e^{4V_1-2V_2}-\Pi |z-\ve p_j|^2\Pi |z-\ve q_j|^2e^{V_2-V_1}
-2|z|^{4N_1}e^{4\tilde{U}_{1,0}-2\tilde{U}_{2,0}}\\
& & +|z|^{2N_1+2N_2}e^{\tilde{U}_{2,0}-\tilde{U}_{1,0}}+ \frac{E_{11}}{\ve^2} \\
& & +f(0,z)e^{2\tilde{U}_{1,0}-\tilde{U}_{2,0}}(2\psi_1-\psi_2)\\
& &+\frac{f(0,z)}{2}e^{2U_{1,0}-U_{2,0}}(2\Psi_{0,1}-\Psi_{0,2})^2
+f_\ve(0,z)e^{2\tilde{U}_{1,0}-\tilde{U}_{2,0}}
(2\Psi_{0,1}-\Psi_{0,2})\\
&&+\frac{f_{\ve\ve}(0,z)}{2}e^{2\tilde{U}_{1,0}-\tilde{U}_{2,0}} ,
\end{eqnarray*}
\begin{eqnarray*}
& \ \ & E_2 = 2\Pi |z-\ve q_j|^4e^{4V_2-6V_1}-3\Pi |z-\ve p_j|^2\Pi |z-\ve q_j|^2e^{V_2-V_1}
-2|z|^{4N_2}e^{4\tilde{U}_{2,0}-6\tilde{U}_{1,0}}\\
& & +3|z|^{2N_1+2N_2}e^{\tilde{U}_{2,0}-\tilde{U}_{1,0}} +\frac{E_{22}}{\ve^2} \\
& &
+g(0,z)e^{2\tilde{U}_{2,0}-3\tilde{U}_{1,0}}(2\psi_2-3\psi_1)\\
&&+\frac{g(0,z)}{2}e^{2U_{2,0}-3U_{1,0}}(2\Psi_{0,2}-3\Psi_{0,1})^2
+g_\ve(0,z)e^{2\tilde{U}_{2,0}-3\tilde{U}_{1,0}}
(2\Psi_{0,2}-3\Psi_{0,1})\\
&&+\frac{g_{\ve\ve}(0,z)}{2}e^{2\tilde{U}_{2,0}-3\tilde{U}_{1,0}}.
\end{eqnarray*}

Here
\begin{eqnarray*}
& & E_{11}= \\
 & &  -f(\ve,z)e^{2V_1-V_2}
+f(0,z)e^{2\tilde{U}_{1,\mathbf{b}}-\tilde{U}_{2,\mathbf{b}}}+\ve f(0,z)e^{2\tilde{U}_{1,0}-\tilde{U}_{2,0}}(2\Psi_{0,1}-\Psi_{0,2})
+\ve f_\ve(0,z)e^{2\tilde{U}_{1,0}-\tilde{U}_{2,0}}  \\
&  &+\sum_{i=3}^{14}\ve[f(0,z)e^{2U_{1,0}-U_{2,0}}(2\Psi_{i,1}-\Psi_{i,2})
+f(0,z)e^{2U_{1,0}-U_{2,0}}(2\Psi_{0,1}-\Psi_{0,2})(2Z_{i,1}-Z_{i,2})]\mathbf{a}_i \\
& &+\sum_{i=3}^{14}\ve f_\ve(0,z)e^{2U_{1,0}-U_{2,0}}(2Z_{i,1}-Z_{i,2})\mathbf{a}_i
\end{eqnarray*}
and
\begin{eqnarray*}
& & E_{22}= \\
&&-g(\ve,z)e^{2V_2-3V_1}
+g(0,z)e^{2\tilde{U}_{2,\mathbf{b}}-3\tilde{U}_{1,\mathbf{b}}}+\ve g(0,z)e^{2\tilde{U}_{2,0}-3\tilde{U}_{1,0}}(2\Psi_{0,2}-3\Psi_{0,1})\\
&&+\ve g_\ve(0,z)e^{2\tilde{U}_{2,0}-3\tilde{U}_{1,0}}\\
&&+\sum_{i=3}^{14}\ve[g(0,z)e^{2U_{2,0}-3U_{1,0}}(2\Psi_{i,2}-3\Psi_{i,1})
+g(0,z)e^{2U_{2,0}-3U_{1,0}}(2\Psi_{0,2}-3\Psi_{0,1})(2Z_{i,2}-3Z_{i,1})]\mathbf{a}_i\\
&&+\sum_{i=3}^{14}\ve g_\ve(0,z)e^{2U_{2,0}-3U_{1,0}}(2Z_{i,2}-3Z_{i,1})\mathbf{a}_i.
\end{eqnarray*}

By Taylor's expansion, we have

\begin{equation}
\label{E1000b}
 \begin{array}{lll}
& & E_1=-f(0,z)e^{2\tilde{U}_{1,0}-\tilde{U}_{2,0}} \Big\{ \frac{1}{\ve}\displaystyle \sum _{i=3}^{14}(2\Psi_{i,1}-\Psi_{i,2})\mathbf{a}_i\displaystyle \sum _{j=3}^{14}(2Z_{j,1}-Z_{j,2})\mathbf{a}_j \\
&&+(2\psi_{1}-\psi_{2})\displaystyle\sum_{i=3}^{14}(2Z_{i,1}-Z_{i,2})\mathbf{a}_i
+\frac{1}{2}(2\Psi_{0,1}-\Psi_{0,2})^2\displaystyle \sum _{i=3}^{14}(2Z_{i,1}-Z_{i,2})\mathbf{a}_i\\
&&+(2\Psi_{0,1}-\Psi_{0,2})\displaystyle\sum_{i=3}^{14}(2\Psi_{i,1}-\Psi_{i,2})\mathbf{a}_i
+(2\Psi_{0,1}-\Psi_{0,2})\displaystyle \sum _{i=3}^{14}(2Z_{i,1}-Z_{i,2})\mathbf{a}_i \Big\}\\
&&-\frac{f(0,z)}{2\ve}\displaystyle \sum _{i,j=3}^{14}\partial^2_{\mathbf{a}_i\mathbf{a}_j}( e^{(2\tilde{U}_{1,0}-\tilde{U}_{2,0})})\left( 2\Psi _{0,1}-\Psi _{0,2}\right)\mathbf{a}_i\mathbf{a}_j\\
&&-f_{\ve}(0,z)\Big\{
 e^{2\tilde{U}_{1,0}-\tilde{U}_{2,0}}(2\Psi_{0,1}-\Psi_{0,2})\displaystyle \sum_{i=3}^{14}(2Z_{i,1}-Z_{i,2})\mathbf{a}_i
+e^{2\tilde{U}_{1,0}-\tilde{U}_{2,0}}\displaystyle\sum_{i=3}^{14}(2\Psi_{i,1}
-\Psi_{i,2})\mathbf{a}_i\\
&&+\frac{1}{2\ve}\displaystyle \sum _{i,j=3}^{14}\partial^2_{\mathbf{a}_i\mathbf{a}_j}( e^{(2\tilde{U}_{1,0}-\tilde{U}_{2,0})})\mathbf{a}_i\mathbf{a}_j\Big\}
-\frac{1}{2}f_{\ve \ve}(0,z)e^{2\tilde{U}_{1,0}-\tilde{U}_{2,0}}\displaystyle \sum_{i=3}^{14}(2Z_{i,1}-Z_{i,2})\mathbf{a}_i\\
&&+4f(0,z)^2 e^{4\tilde{U}_{1,0}-2\tilde{U}_{2,0}}\displaystyle \sum_{i=3}^{14}(2Z_{i,1}-Z_{i,2})\mathbf{a}_i
-f(0,z)g(0,z)\displaystyle \sum_{i=3}^{14} e^{\tilde{U}_{2,0}-\tilde{U}_{1,0}}(Z_{i,2}-Z_{i,1})\mathbf{a}_i\\
&&+O(\ve)+O(\ve ^2+|\mathbf{a}|^2),
\end{array}
\end{equation}
and
\begin{equation}
\label{E2000b}
 \begin{array}{lll}
& & E_2=-g(0,z)e^{2\tilde{U}_{2,0}-3\tilde{U}_{1,0}} \Big\{ \frac{1}{\ve}\displaystyle \sum _{i=3}^{14}(2\Psi_{i,2}-3\Psi_{i,1})\mathbf{a}_i\displaystyle \sum _{j=3}^{14}(2Z_{j,2}-3Z_{j,1})\mathbf{a}_j \\
&&+(2\psi_{2}-3\psi_{1})\displaystyle\sum_{i=3}^{14}(2Z_{i,2}-3Z_{i,1})\mathbf{a}_i
+\frac{1}{2}(2\Psi_{0,2}-3\Psi_{0,1})^2\displaystyle \sum _{i=3}^{14}(2Z_{i,2}-3Z_{i,1})\mathbf{a}_i\\
&&+(2\Psi_{0,2}-3\Psi_{0,1})\displaystyle\sum_{i=3}^{14}(2\Psi_{i,2}-3\Psi_{i,1})\mathbf{a}_i
+(2\Psi_{0,2}-3\Psi_{0,1})\displaystyle \sum _{i=3}^{14}(2Z_{i,2}-3Z_{i,1})\mathbf{a}_i \Big\}\\
&&-\frac{g(0,z)}{2\ve}\displaystyle \sum _{i,j=3}^{14}\partial^2_{\mathbf{a}_i\mathbf{a}_j}( e^{(2\tilde{U}_{2,0}-3\tilde{U}_{1,0})})( 2\Psi _{0,2}-3\Psi _{0,1})\mathbf{a}_i\mathbf{a}_j\\
&&-g_{\ve}(0,z)\Big\{
 e^{2\tilde{U}_{2,0}-3\tilde{U}_{1,0}}(2\Psi_{0,2}-3\Psi_{0,1})\displaystyle \sum_{i=3}^{14}(2Z_{i,2}-3Z_{i,1})\mathbf{a}_i
+e^{2\tilde{U}_{2,0}-3\tilde{U}_{1,0}}\displaystyle\sum_{i=3}^{14}(2\Psi_{i,2}
-3\Psi_{i,1})\mathbf{a}_i\\
&&+\frac{1}{2\ve}\displaystyle \sum _{i,j=3}^{14}\partial^2_{\mathbf{a}_i\mathbf{a}_j}( e^{(2\tilde{U}_{2,0}-3\tilde{U}_{1,0})})\mathbf{a}_i\mathbf{a}_j\Big\}
-\frac{1}{2}g_{\ve \ve}(0,z)e^{2\tilde{U}_{2,0}-3\tilde{U}_{1,0}}\displaystyle \sum_{i=3}^{14}(2Z_{i,2}-3Z_{i,1})\mathbf{a}_i\\
&&+4g(0,z)^2 e^{4\tilde{U}_{2,0}-6\tilde{U}_{1,0}}\displaystyle \sum_{i=3}^{14}(2Z_{i,2}-3Z_{i,1})\mathbf{a}_i
-3f(0,z)g(0,z)\displaystyle \sum_{i=3}^{14} e^{\tilde{U}_{2,0}-\tilde{U}_{1,0}}(Z_{i,2}-Z_{i,1})\mathbf{a}_i\\
&&+O(\ve)+O(\ve ^2+|\mathbf{a}|^2),
\end{array}
\end{equation}
where $O(\ve)$ denotes all items only involving with $\ve$,  and not with $\mathbf{a}$.

\section{A Nonlinear Projected Problem}\label{sec4}
Similar to Proposition 2.1 in \cite{ALW}, we have the following result:
\begin{proposition}\label{propb}
For $\mathbf{a}$ satisfying (\ref{parab}), there exists a solution $(v,\{m_i\})$ to the following system
\begin{equation}\label{projectb}
\left\{\begin{array}{c}
\Delta v_1+|z|^{2N_1}e^{2\tilde{U}_{1,0}-\tilde{U}_{2,0}}(2v_1-v_2)=G_1+\sum_{i=3}^{14}m_i(v) Z_{i,1}^*\\
\Delta v_2+|z|^{2N_2}e^{2\tilde{U}_{2,0}-3\tilde{U}_{1,0}}(2v_2-3v_1)=G_2+\sum_{i=3}^{14}m_i(v) Z_{i,2}^*\\
\langle\Delta Z_i, v\rangle=0, \mbox{ for }i=1,\cdots,14,
\end{array}
\right.
\end{equation}
where $G=\vect{G_1}{G_2} $ and $m_i(v)$ can be determined by
\begin{equation}\label{e402b}
\langle G+\displaystyle \sum _{i=3}^{14} m_i(v)Z_i^*, Z_j^*\rangle =0, \; \mbox{for} \; j=3,\ldots,14.
\end{equation}
\noindent
Furthermore, $v$ satisfies the following estimate
\begin{equation}
\|v\|_*\leq C\ve,
\end{equation}
for some constant $C$ independent of $\ve$.
\end{proposition}

\bigskip

By Proposition \ref{propb}, the full solvability for \eqref{3b} is reduced to $m_i=0$ for $i=3,\cdots,14$. Since by (\ref{det2b}), $\rm det(\langle Z_i^*,Z_j^*\rangle_{i,j=3,\cdots,14})\neq 0$, and recall the definition of $m_i$ in (\ref{e402b}), $m_i=0$ is equivalent to
 \begin{equation}\label{errorb}
 \int_0^{+\infty}\int_0^{2\pi} G\cdot Z_i^* rd\theta dr=0 \  \mbox{ for }i=3,\cdots,14.
 \end{equation}

To solve (\ref{errorb}), we observe all $N_{ij}$ terms in (\ref{eq4}) and (\ref{eq5}) are small. In fact, we have the following lemma:

\begin{lemma}\label{lemma401}
Let $\vect{v_1}{v_2}$ be a solution of (\ref{projectb}). Then we have the following estimates:
\begin{equation}\label{N(v)b}
\int_{\R^2} (N_{11}(v)+N_{12}(v))Z_{i,1}^*+(N_{21}(v)+N_{22}(v))Z_{i,2}^* dx=O(\ve^2),
\end{equation}
for $i=3,\cdots,14$, where $O(\ve^2)\leq C_1\ve^2$ for some positive $C_1$ independent of ${\bf a}$ provided $|{\bf a}|\leq 1$.
\end{lemma}
\noindent
{\bf Proof:}

The proof is similar to that of Lemma 2.4 in \cite{ALW}.

\qed

From the Taylor expansions in (\ref{E1000b}) and (\ref{E2000b}), we obtain that the error projection can be expressed as
\begin{equation}
(\langle E, Z_{{\bf a}_i}^*\rangle)_{i=3,\cdots,14}=\frac{1}{\ve}\tilde{\bf A}{\bf a}\cdot{\bf a}+\tilde{\bf Q}{\bf a}+O(|{\bf a}|^2)+O(\ve).
\end{equation}

\section{Proof of  Theorem \ref{thm101} under Assumption (i)}\label{sec3.5}

Suppose the {\bf Assumption (i)} holds. By a translation, we might assmue that $\sum_{i=1}^{N_1}p_i=\sum_{j=1}^{N_2}q_j=0$ and $N_1=N_2$ , and we choose $(\xi_1,\xi_2)=(0,0)$ in this section. This case is the reminiscent of $SU(2)$ case, even though, the proof is considerably harder since there are fourteen dimensional kernels instead of a three-dimensional one for the $SU(2)$ case.

\begin{lemma}\label{lemma402b}
Let $\vect{v_1}{v_2}$ be a solution of (\ref{projectb}). The following estimates hold:
\begin{eqnarray}
&&(\langle E, Z_{c_{43,1}}^*\rangle,\langle E, Z_{c_{43,2}}^*\rangle,\langle E, Z_{c_{52,1}}^*\rangle,\langle E, Z_{c_{52,2}}^*\rangle,\\
&& \, \, \, \, \, \, \, \, \, \, \langle E, Z_{c_{53,1}}^*\rangle,\langle E, Z_{c_{53,2}}^*\rangle,\langle E, Z_{c_{54,1}}^*\rangle,\langle E, Z_{c_{54,2}}^*\rangle,\nonumber\\
&& \, \, \, \, \, \, \, \, \, \,\langle E, Z_{c_{61,1}}^*\rangle,\langle E, Z_{c_{61,2}}^*\rangle,
\langle E, Z_{c_{62,1}}^*\rangle,\langle E, Z_{c_{62,2}}^*\rangle)^t\nonumber\\
&&=\tilde{\mathcal{T}}(\mathbf{a})+O(|\mathbf{a}|^2)+O(\ve),\nonumber
\end{eqnarray}
where $\tilde{\mathcal{T}} $ is an $12\times 12$ matrix defined in (\ref{Tb}). Moreover, $\tilde{\mathcal{T}}$ is non-degenerate if $N>0$.
\end{lemma}
\noindent
{\bf Proof:}

 Without loss of generality , we may assume that $ \sum_{j=1}^{N_1} p_j= \sum_{j=1}^{N_2} q_j=0$ and $ N_1=N_2=N$, and denote by $\mu=N+1$. Now we choose the parameters $(\lambda_4,\lambda_5,\xi_1,\xi_2)=(\frac{1}{3\times 2^{10}\mu^6},\frac{1}{15\times 2^{9}\mu^6},0,0)$,  so that we have
 \begin{equation}
 e^{\tilde{U}_{1,0}}=\rho_{1,G}=\frac{1}{2}\rho_1=\frac{45\times 2^{9} \mu^6}{(1+r^{2\mu})^6},\ \  e^{\tilde{U}_{2,0}}=\rho_{2,G}=\frac{1}{4}\rho_2=\frac{3^3\times 5^2\times 2^{15} \mu^{10}}{(1+r^{2\mu})^{10}}.
 \end{equation}
Since $\sum_{i=1}^{N}p_i=\sum_{i=1}^{N}q_i=0$, we have $f_{\ve}(0,z)=0$, $\Psi _{0,1}=\Psi_{0,2}=0$. By \eqref{Psi1b}, we have $\Psi _{i,1}=\Psi _{i,2}=0$, $i=3, \ldots,14$. Recall that
\begin{equation}
E=\vect{E_1}{E_2} ,
\end{equation}
where by \eqref{E1000b}, and \eqref{E2000b}, we have
\begin{eqnarray*}
E_1&=&-|z|^{2N}\frac{\rho_{1,G}^2}{\rho_{2,G}} (2\psi_1-\psi_2)\displaystyle \sum_{i=3}^{14}(2Z_{i,1}-Z_{i,2})\mathbf{a}_i-\frac{1}{2} f_{\ve \ve}(0,z)\frac{\rho_{1,G}^2}{\rho_{2,G}}\displaystyle \sum_{i=3}^{14}(2Z_{i,1}-Z_{i,2})\mathbf{a}_i\\
&+&4|z|^{4N}\frac{\rho_{1,G}^4}{\rho_{2,G}^2}\displaystyle \sum _{i=3}^{14}(2Z_{i,1}-Z_{i,2})\mathbf{a}_i-|z|^{4N}\frac{\rho_{2,G}}{\rho_{1,G}} \displaystyle \sum_{i=3}^{14}(Z_{i,2}-Z_{i,1})\mathbf{a}_i\\
&+&O(\ve)+O(\ve ^2+|\mathbf{a}|^2),
\end{eqnarray*}
\begin{eqnarray*}
E_2&=&-|z|^{2N}\frac{\rho_{2,G}^2}{\rho_{1,G}^3} (2\psi_2-3\psi_1)\displaystyle \sum _{i=3}^{14}(2Z_{i,2}-3Z_{i,1})\mathbf{a}_i-\frac{1}{2} g_{\ve \ve}(0,z)\frac{\rho_{2,G}^2}{\rho_{1,G}^3}\displaystyle \sum_{i=3}^{14}(2Z_{i,2}-3Z_{i,1})\mathbf{a}_i\\
&+&4|z|^{4N}\frac{\rho_{2,G}^4}{\rho_{1,G}^6}\displaystyle \sum _{i=3}^{14}(2Z_{i,2}-3Z_{i,1})\mathbf{a}_i-3|z|^{4N}\frac{\rho_{2,G}}{\rho_{1,G}} \displaystyle \sum_{i=3}^{14}(Z_{i,2}-Z_{i,1})\mathbf{a}_i\\
&+&O(\ve)+O(\ve ^2+|\mathbf{a}|^2),
\end{eqnarray*}
where $f_{\ve \ve} (0,z)$ is

\begin{eqnarray}\label{fveveb}
&&f_{\ve\ve}(0,z)=2|z|^{2(N-1)}(\sum_i|p_i|^2+2\sum_{i\neq j}(p_{i1}\cos\theta+p_{i2}\sin\theta)(p_{j1}\cos\theta+p_{j2}\sin\theta))\nonumber\\
&&=2|z|^{2(N-1)}(|\sum_i p_i|^2+\sum_{i\neq j}(p_{i1}p_{j1}-p_{i2}p_{j2})\cos {2\theta}+(p_{i1}p_{j2}+p_{i2}p_{j1})\sin{2\theta})\nonumber\\
&&=2|z|^{2(N-1)}(\sum_{i\neq j}(p_{i1}p_{j1}-p_{i2}p_{j2})\cos {2\theta}+(p_{i1}p_{j2}+p_{i2}p_{j1})\sin{2\theta}).
\end{eqnarray}
Similarly we can get the expression for $g_{\ve\ve}(0,z)$.

Since
\begin{eqnarray*}
& & \int h(r)\cos{2\theta}(2Z_{i,1}-Z_{i,2})(2Z_{j,1}-Z_{j,2})rdrd\theta\\
&= & \int h(r)\sin{2\theta}(2Z_{i,1}-Z_{i,2})(2Z_{j,1}-Z_{j,2})rdrd\theta=0,
\end{eqnarray*}
and
\begin{eqnarray*}
& & \int h(r)\cos{2\theta}(Z_{i,2}-Z_{i,1})(Z_{j,2}-Z_{j,1})rdrd\theta\\
&= & \int h(r)\sin{2\theta}(Z_{i,2}-Z_{i,1})(Z_{j,2}-Z_{j,1})rdrd\theta=0,
\end{eqnarray*}
for $i,j=3,\cdots,14$, from (\ref{fveveb}), we have
\begin{equation}\label{fintb}
\int_0^\infty\int_0^{2\pi}f_{\ve\ve}(0,z)e^{2\tilde{U}_{1,0}-\tilde{U}_{2,0}}
(2Z_{i,1}-Z_{i,2})(2Z_{j,1}-Z_{j,2})rdrd\theta=0,
\end{equation}
and
\begin{equation}\label{gintb}
\int_0^\infty\int_0^{2\pi}g_{\ve\ve}(0,z)e^{2\tilde{U}_{2,0}-3\tilde{U}_{1,0}}
(2Z_{i,2}-3Z_{i,1})(2Z_{j,2}-3Z_{j,1})rdrd\theta=0,
\end{equation}
for $i,j=3,\cdots,14$. Note that $f_{\ve\ve}(0,z)=g_{\ve \ve}(0,z)=0$ if $N=1$.

Another important observation is the following:
\begin{eqnarray}\label{psi1psi10b}
&&\int_0^\infty\int_0^{2\pi}r^{2N}e^{2\tilde{U}_{1,0}-\tilde{U}_{2,0}}
(2\psi_1-\psi_2)(2Z_{i,1}-Z_{i,2})(2Z_{j,1}-Z_{j,2})rdrd\theta\\
&&=\int_0^\infty\int_0^{2\pi}f(0,z)e^{2\tilde{U}_{1,0}-\tilde{U}_{2,0}}
(2\psi_{1}^0-\psi_{2}^0)(2Z_{i,1}-Z_{i,2})(2Z_{j,1}-Z_{j,2})rdrd\theta,\nonumber\\
&&\nonumber\\
&&\int_0^\infty\int_0^{2\pi}r^{2N}e^{2\tilde{U}_{2,0}-3\tilde{U}_{1,0}}
(2\psi_2-3\psi_1)(2Z_{i,2}-3Z_{i,1})(2Z_{j,2}-3Z_{j,1})rdrd\theta \label{psi2psi20b}
\\
&&=\int_0^\infty\int_0^{2\pi}g(0,z)e^{2\tilde{U}_{2,0}-3\tilde{U}_{1,0}}
(2\psi_{2}^0-3\psi_{1}^0)(2Z_{i,2}-3Z_{i,1})(2Z_{j,2}-3Z_{j,1})rdrd\theta,\nonumber
\end{eqnarray}
where $(\psi_{1}^0,\psi_{2}^0)$ is the radial solution of the following system:
\begin{equation}
\left\{\begin{array}{c}
\Delta \psi_1+|z|^{2N_1}e^{2\tilde{U}_{1,0}-\tilde{U}_{2,0}}(2\psi_1-\psi_2)=
2|z|^{4N_1}e^{4\tilde{U}_{1,0}-2\tilde{U}_{2,0}}
-|z|^{2(N_1+N_2)}e^{\tilde{U}_{2,0}-\tilde{U}_{1,0}}\\
\Delta \psi_2+|z|^{2N_2}e^{2\tilde{U}_{2,0}-3\tilde{U}_{1,0}}(2\psi_2-3\psi_1)=
2|z|^{4N_2}e^{4\tilde{U}_{2,0}-6\tilde{U}_{1,0}}
-3|z|^{2(N_1+N_2)}e^{\tilde{U}_{2,0}-\tilde{U}_{1,0}}.\\
\end{array}
\right.
\end{equation}
Obviously, $(\psi_1^0,\psi_2^0)$ is the radial part of $(\psi_{0,1},\psi_{0,2})$. Because of this observation, when dealing with the $O(\ve^2)$ approximation, we only need to consider the radial part of the solutions.

In fact we can choose  $\psi_{1}^0=\psi,\ \ \psi_{2}^0=\frac{5}{3}\psi$ such that  $\psi$ is the solution of the following ODE:
\begin{equation}
\Delta \psi+\frac{8(N+1)^2r^{2N}}{(1+r^{2N+2})^2}\psi=r^{4N}\frac{3\times 2^6 (N+1)^4}{(1+r^{2N+2})^4}.
\end{equation}

Combining (\ref{fintb}), (\ref{gintb}), (\ref{psi1psi10b}) and (\ref{psi2psi20b}), one can get the following:
\begin{eqnarray*}
&&\int E\cdot Z_{k}^* rdrd\theta\\
&=&\int_0^\infty\int_0^{2\pi}\sum_{i=3}^{14}\Big(\Big[2|z|^{4N_1}e^{4\tilde{U}_{1,0}-2\tilde{U}_{2,0}}(4Z_{i,1}-2Z_{i,2})\mathbf{a}_i
-|z|^{2N_1+2N_2}e^{\tilde{U}_{2,0}-\tilde{U}_{1,0}}(Z_{i,2}-Z_{1,i})\mathbf{a}_i\\
&-&f(0,z)e^{2\tilde{U}_{1,0}-\tilde{U}_{2,0}}(2\psi_{0,1}-\psi_{0,2})(2Z_{i,1}-Z_{i,2})\mathbf{a}_i\Big]Z_{k,1}^*\\
&+&\Big[4|z|^{4N_2}e^{4\tilde{U}_{2,0}-6\tilde{U}_{1,0}}(2Z_{i,2}-3Z_{i,1})\mathbf{a}_i
-3|z|^{2N_1+2N_2}e^{\tilde{U}_{2,0}-\tilde{U}_{1,0}}(Z_{i,2}-Z_{i,1})\mathbf{a}_i\\
&-&g(0,z)e^{2\tilde{U}_{2,0}-3\tilde{U}_{1,0}}(2\psi_{0,2}-3\psi_{0,1})(2Z_{i,2}-3Z_{i,1})\mathbf{a}_i\Big]Z_{k,2}^*\Big)rdrd\theta\\
&+&O(|\mathbf{a}|^2+|\ve |^2)+O(\ve),
\end{eqnarray*}
where $O(\ve)$ is independent of $\mathbf{a}$.

\medskip

Thus we can get that for $i=1,2$
\begin{eqnarray*}
&&\int E\cdot Z_{c_{43,i}}^* rdrd\theta\\
&=&\int \Big(4r^{4N}\frac{\rho_{1,G}^4}{\rho_{2,G}^2}\Big[(2Z_{c_{43,i},1}-Z_{c_{43,i},2})^2c_{43,i}+(2Z_{c_{43,i},1}-Z_{c_{43,i},2})
(2Z_{c_{54,i},1}-Z_{c_{54,i},2})c_{54,i}\Big]\\
&-&r^{4N}\frac{\rho_{2,G}}{\rho_{1,G}}\Big[(Z_{c_{43,i},2}-Z_{c_{43,i},1})(2Z_{c_{43,i},1}-Z_{c_{43,i},2})c_{43,i}
+(Z_{c_{54,i},2}-Z_{c_{54,i},1})(2Z_{c_{43,i},1}-Z_{c_{43,i},2})c_{54,i}\Big]\\
&-&r^{2N}\frac{\rho_{1,G}^2}{\rho_{2,G}}(2\psi_{1,0}-\psi_{2,0})\Big[(2Z_{c_{43,i},1}-Z_{c_{43,i},2})^2c_{43,i}+
(2Z_{c_{43,i},1}-Z_{c_{43,i},2})(2Z_{c_{54,i},1}-Z_{c_{54,i},2})c_{54,i}\Big]\\
&+&4r^{4N}\frac{\rho_{2,G}^4}{\rho_{1,G}^6}\Big[\frac{1}{3}(2Z_{c_{43,i},2}-3Z_{c_{43,i},1})^2c_{43,i}+
\frac{1}{3}(2Z_{c_{43,i},2}-3Z_{c_{43,i},1})(2Z_{c_{54,i},2}-3Z_{c_{54,i},1})c_{54,i}\Big]\\
&-&r^{4N}\frac{\rho_{2,G}}{\rho_{1,G}}\Big[(Z_{c_{43,i},2}-Z_{c_{43,i},1})(2Z_{c_{43,i},2}-3Z_{c_{43,i},1})c_{43,i}
+(Z_{c_{54,i},2}-Z_{c_{54,i},1})(2Z_{c_{43,i},2}-3Z_{c_{43,i},1})c_{54,i}\Big]\\
&-&r^{2N}\frac{\rho_{2,G}^2}{\rho_{1,G}^3}(2\psi_{2,0}-3\psi_{1,0})
\Big[\frac{1}{3}(2Z_{c_{43,i},2}-3Z_{c_{43,i},1})^2c_{43,i}\\
&+&\frac{1}{3}(2Z_{c_{43,i},2}-3Z_{c_{43,i},1})(2Z_{c_{54,i},2}-3Z_{c_{54,i},1})c_{54,i}\Big]
\Big)rdrd\theta\\
&=&\int \Big[\Big(4r^{4N}\frac{\rho_{1,G}^4}{\rho_{2,G}^2}(2Z_{c_{43,i},1}-Z_{c_{43,i},2})^2
-r^{4N}\frac{\rho_{2,G}}{\rho_{1,G}}(Z_{c_{43,i},2}-Z_{c_{43,i},1})(2Z_{c_{43,i},1}-Z_{c_{43,i},2})\\
&+&\frac{4}{3}r^{4N}\frac{\rho_{2,G}^4}{\rho_{1,G}^6}(2Z_{c_{43,i},2}-3Z_{c_{43,i},1})^2
-r^{4N}\frac{\rho_{2,G}}{\rho_{1,G}}(Z_{c_{43,i},2}-Z_{c_{43,i},1})(2Z_{c_{43,i},2}-3Z_{c_{43,i},1})\Big)\\
&-&\frac{1}{3}r^{2N}\Big(\frac{\rho_{1,G}^2}{\rho_{2,G}}(2Z_{c_{43,i},1}-Z_{c_{43,i},2})^2+\frac{1}{3}\frac{\rho_{2,G}^2}{\rho_{1,G}^3}(2Z_{c_{43,i},2}-3Z_{c_{43,i},1})^2\Big)
\psi\Big)c_{43,i}\\
&+&\Big(\Big(4r^{4N}\frac{\rho_{1,G}^4}{\rho_{2,G}^2}(2Z_{c_{43,i},1}-Z_{c_{43,i},2})(2Z_{c_{54,i},1}-Z_{c_{54,i},2})
-r^{4N}\frac{\rho_{2,G}}{\rho_{1,G}}(Z_{c_{54,i},2}-Z_{c_{54,i},1})(2Z_{c_{43,i},1}-Z_{c_{43,i},2})\\
&+&\frac{4r^{4N}}{3}\frac{\rho_{2,G}^4}{\rho_{1,G}^6}(2Z_{c_{43,i},2}-3Z_{c_{43,i},1})
(2Z_{c_{54,i},2}-3Z_{c_{54,i},1})
-\frac{\rho_{2,G}r^{4N}}{\rho_{1,G}}(Z_{c_{54,i},2}-Z_{c_{54,i},1})(2Z_{c_{43,i},2}-3Z_{c_{43,i},1})\Big)\\
&-&\frac{r^{2N}}{3}\Big(\frac{\rho_{1,G}^2}{\rho_{2,G}} (2Z_{c_{43,i},1}-Z_{c_{43,i},2})(2Z_{c_{54,i},1}-Z_{c_{54,i},2})\\
&+&\frac{1}{3}\frac{\rho_{2,G}^2}{\rho_{1,G}^3}(2Z_{c_{43,i},2}-3Z_{c_{43,i},1})(2Z_{c_{54,i},2}-3Z_{c_{54,i},1})\Big)
\psi \Big)c_{54,i}\Big]rdrd\theta\\
&=&\pi(J_1+\int q_1\psi rdr)c_{43,i}+(J_2+\int q_7\psi rdr)c_{54,i}\\
&+&O(\ve)+O(|\mathbf{a}|^2+\ve ^2),
\end{eqnarray*}
similarly, we can get that
\begin{eqnarray*}
&&\int E\cdot Z_{c_{54,i}}^* rdrd\theta\\
&=&\int \Big(4r^{4N}\frac{\rho_{1,G}^4}{\rho_{2,G}^2}\Big[(2Z_{c_{54,i},1}-Z_{c_{54,i},2})^2c_{54,i}+(2Z_{c_{43,i},1}-Z_{c_{43,i},2})
(2Z_{c_{54,i},1}-Z_{c_{54,i},2})c_{43,i}\Big]\\
&-&r^{4N}\frac{\rho_{2,G}}{\rho_{1,G}}\Big[(Z_{c_{54,i},2}-Z_{c_{54,i},1})(2Z_{c_{54,i},1}-Z_{c_{54,i},2})c_{54,i}
+(Z_{c_{43,i},2}-Z_{c_{43,i},1})(2Z_{c_{54,i},1}-Z_{c_{54,i},2})c_{43,i}\Big]\\
&-&r^{2N}\frac{\rho_{1,G}^2}{\rho_{2,G}}(2\psi_{1,0}-\psi_{2,0})\Big[(2Z_{c_{54,i},1}-Z_{c_{54,i},2})^2c_{54,i}+
(2Z_{c_{43,i},1}-Z_{c_{43,i},2})(2Z_{c_{54,i},1}-Z_{c_{54,i},2})c_{43,i}\Big]\\
&+&4r^{4N}\frac{\rho_{2,G}^4}{\rho_{1,G}^6}\Big[\frac{1}{3}(2Z_{c_{54,i},2}-3Z_{c_{54,i},1})^2c_{54,i}+
\frac{1}{3}(2Z_{c_{54,i},2}-3Z_{c_{54,i},1})(2Z_{c_{43,i},2}-3Z_{c_{43,i},1})c_{43,i}\Big]\\
&-&r^{4N}\frac{\rho_{2,G}}{\rho_{1,G}}\Big[(Z_{c_{54,i},2}-Z_{c_{54,i},1})(2Z_{c_{54,i},2}-3Z_{c_{54,i},1})c_{54,i}
+(Z_{c_{43,i},2}-Z_{c_{43,i},1})(2Z_{c_{54,i},2}-3Z_{c_{54,i},1})c_{43,i}\Big]\\
&-&r^{2N}\frac{\rho_{2,G}^2}{\rho_{1,G}^3}(2\psi_{2,0}-3\psi_{1,0})\Big[\frac{1}{3}(2Z_{c_{54,i},2}-3Z_{c_{54,i},1})^2c_{54,i}\\
&+&\frac{1}{3} (2Z_{c_{43,i},2}-3Z_{c_{43,i},1})(2Z_{c_{54,i},2}-3Z_{c_{54,i},1})c_{43,i}\Big]\Big)rdrd\theta\\
&=&\int \Big[\Big(4r^{4N}\frac{\rho_{1,G}^4}{\rho_{2,G}^2}(2Z_{c_{54,i},1}-Z_{c_{54,i},2})^2
-r^{4N}\frac{\rho_{2,G}}{\rho_{1,G}}(Z_{c_{54,i},2}-Z_{c_{54,i},1})(2Z_{c_{54,i},1}-Z_{c_{54,i},2})\\
&+&\frac{4}{3} r^{4N}\frac{\rho_{2,G}^4}{\rho_{1,G}^6}(2Z_{c_{54,i},2}-3Z_{c_{54,i},1})^2
-r^{4N}\frac{\rho_{2,G}}{\rho_{1,G}}(Z_{c_{54,i},2}-Z_{c_{54,i},1})(2Z_{c_{54,i},2}-3Z_{c_{54,i},1})\Big)\\
&-&\frac{1}{3}r^{2N}\Big(\frac{\rho_{1,G}^2}{\rho_{2,G}}(2Z_{c_{54,i},1}-Z_{c_{54,i},2})^2+\frac{1}{3}
\frac{\rho_{2,G}^2}{\rho_{1,G}^3}(2Z_{c_{54,i},2}-3Z_{c_{54,i},1})^2\Big)\psi
\Big)c_{54,i}\\
&+&\Big(4r^{4N}\frac{\rho_{1,G}^4}{\rho_{2,G}^2}(2Z_{c_{43,i},1}-Z_{c_{43,i},2})(2Z_{c_{54,i},1}-Z_{c_{54,i},2})-r^{4N}\frac{\rho_{2,G}}{\rho_{1,G}}(Z_{c_{54,i},2}-Z_{c_{54,i},1})(2Z_{c_{43,i},1}-Z_{c_{43,i},2})\\
&+&\frac{4r^{4N}}{3}\frac{\rho_{2,G}^4}{\rho_{1,G}^6}(2Z_{c_{43,i},2}-3Z_{c_{43,i},1})
(2Z_{c_{54,i},2}-3Z_{c_{54,i},1})
-\frac{\rho_{2,G}r^{4N}}{\rho_{1,G}}(Z_{c_{43,i},2}-Z_{c_{43,i},1})(2Z_{c_{54,i},2}-3Z_{c_{54,i},1})\\
&-&\frac{1}{3}r^{2N}\Big(\frac{\rho_{1,G}^2}{\rho_{2,G}}(2Z_{c_{54,i},1}-Z_{c_{54,i},2})
(2Z_{c_{43,i},1}-Z_{c_{43,i},2})\\
&+&\frac{1}{3}\frac{\rho_{2,G}^2}{\rho_{1,G}^3}(2Z_{c_{54,i},2}-3Z_{c_{54,i},1})(2Z_{c_{43,i},2}-3Z_{c_{43,i},1})\Big)
\psi \Big)c_{43,i}\Big]rdrd\theta\\
&=&\pi(J_3+\int q_4 \psi r dr)c_{43,i}+(J_4+\int q_7 \psi r dr)c_{54,i}\\
&+&O(\ve)+O(|\mathbf{a}|^2+\ve ^2),
\end{eqnarray*}
and

\begin{eqnarray*}
&&\int E\cdot Z_{c_{52,i}}^* rdrd\theta\\
&=&\int \Big(4r^{4N}\frac{\rho_{1,G}^4}{\rho_{2,G}^2}\Big[(2Z_{c_{52,i},1}-Z_{c_{52,i},2})^2c_{52,i}\Big]
-r^{4N}\frac{\rho_{2,G}}{\rho_{1,G}}\Big[(Z_{c_{52,i},2}-Z_{c_{52,i},1})(2Z_{c_{52,i},1}-Z_{c_{52,i},2})c_{52,i}
\Big]\\
&-&r^{2N}\frac{\rho_{1,G}^2}{\rho_{2,G}}(2\psi_{1,0}-\psi_{2,0})
\Big[(2Z_{c_{52,i},1}-Z_{c_{52,i},2})^2c_{52,i}\Big]\\
&+&4r^{4N}\frac{\rho_{2,G}^4}{\rho_{1,G}^6}\Big[\frac{1}{3}(2Z_{c_{52,i},2}-3Z_{c_{52,i},1})^2c_{52,i}\Big]
-r^{4N}\frac{\rho_{2,G}}{\rho_{1,G}}\Big[(Z_{c_{52,i},2}-Z_{c_{52,i},1})(2Z_{c_{52,i},2}-3Z_{c_{52,i},1})c_{52,i}
\Big]\\
&-&r^{2N}\frac{\rho_{2,G}^2}{\rho_{1,G}^3}(2\psi_{2,0}-3\psi_{1,0})\Big[\frac{1}{3}
(2Z_{c_{52,i},2}-3Z_{c_{52,i},1})^2c_{52,i}
\Big]\Big)rdrd\theta\\
&=&\int \Big[\Big(4r^{4N}\frac{\rho_{1,G}^4}{\rho_{2,G}^2}(2Z_{c_{52,i},1}-Z_{c_{52,i},2})^2
-r^{4N}\frac{\rho_{2,G}}{\rho_{1,G}}(Z_{c_{52,i},2}-Z_{c_{52,i},1})(2Z_{c_{52,i},1}-Z_{c_{52,i},2})\\
&+&\frac{4}{3}r^{4N}\frac{\rho_{2,G}^4}{\rho_{1,G}^6}(2Z_{c_{52,i},2}-3Z_{c_{52,i},1})^2
-r^{4N}\frac{\rho_{2,G}}{\rho_{1,G}}(Z_{c_{52,i},2}-Z_{c_{52,i},1})(2Z_{c_{52,i},2}-3Z_{c_{52,i},1})\Big)\\
&-&\frac{1}{3}r^{2N}\Big(\frac{\rho_{1,G}^2}{\rho_{2,G}}(2Z_{c_{52,i},1}-Z_{c_{52,i},2})^2
+\frac{1}{3}\frac{\rho_{2,G}^2}{\rho_{1,G}^3}(2Z_{c_{52,i},2}-3Z_{c_{52,i},1})^2\Big)\psi
\Big)c_{52,i}\Big]r drd\theta\\
&=&\pi(J_5+\int q_2 \psi r dr)c_{52,i}+O(\ve)+O(|\mathbf{a}|^2+\ve ^2),
\end{eqnarray*}

\begin{eqnarray*}
&&\int E\cdot Z_{c_{53,i}}^* rdrd\theta\\
&=&\int \Big(4r^{4N}\frac{\rho_{1,G}^4}{\rho_{2,G}^2}\Big[(2Z_{c_{53,i},1}-Z_{c_{53,i},2})^2c_{53,i}\Big]
-r^{4N}\frac{\rho_{2,G}}{\rho_{1,G}}\Big[(Z_{c_{53,i},2}-Z_{c_{53,i},1})(2Z_{c_{53,i},1}-Z_{c_{53,i},2})c_{53,i}
\Big]\\
&-&r^{2N}\frac{\rho_{1,G}^2}{\rho_{2,G}}(2\psi_{1,0}-\psi_{2,0})
\Big[(2Z_{c_{53,i},1}-Z_{c_{53,i},2})^2c_{53,i}\Big]\\
&+&4r^{4N}\frac{\rho_{2,G}^4}{\rho_{1,G}^6}\Big[\frac{1}{3}(2Z_{c_{53,i},2}-3Z_{c_{53,i},1})^2c_{53,i}\Big]
-r^{4N}\frac{\rho_{2,G}}{\rho_{1,G}}\Big[(Z_{c_{53,i},2}-Z_{c_{53,i},1})(2Z_{c_{53,i},2}-3Z_{c_{53,i},1})c_{53,i}
\Big]\\
&-&r^{2N}\frac{\rho_{2,G}^2}{\rho_{1,G}^3}(2\psi_{2,0}-3\psi_{1,0})\Big[\frac{1}{3}
(2Z_{c_{53,i},2}-3Z_{c_{53,i},1})^2c_{53,i}\Big]
\Big)rdrd\theta\\
&=&\int \Big[\Big(4r^{4N}\frac{\rho_{1,G}^4}{\rho_{2,G}^2}(2Z_{c_{53,i},1}-Z_{c_{53,i},2})^2
-r^{4N}\frac{\rho_{2,G}}{\rho_{1,G}}(Z_{c_{53,i},2}-Z_{c_{53,i},1})(2Z_{c_{53,i},1}-Z_{c_{53,i},2})\\
&+&\frac{4}{3}r^{4N}\frac{\rho_{2,G}^4}{\rho_{1,G}^6}(2Z_{c_{53,i},2}-3Z_{c_{53,i},1})^2
-r^{4N}\frac{\rho_{2,G}}{\rho_{1,G}}(Z_{c_{53,i},2}-Z_{c_{53,i},1})(2Z_{c_{53,i},2}-3Z_{c_{53,i},1})\Big)\\
&-&\frac{1}{3}r^{2N}\Big(\frac{\rho_{1,G}^2}{\rho_{2,G}}(2Z_{c_{53,i},1}-Z_{c_{53,i},2})^2
+\frac{1}{3}\frac{\rho_{2,G}^2}{\rho_{1,G}^3}(2Z_{c_{53,i},2}-3Z_{c_{53,i},1})^2\Big)\psi
\Big)c_{53,i}\\
&=&\pi(J_{6}+\int q_3 \psi r dr)c_{53,i}+O(\ve)+O(|\mathbf{a}|^2+\ve ^2),
\end{eqnarray*}
\begin{eqnarray*}
&&\int E\cdot Z_{c_{61,i}}^* rdrd\theta\\
&=&\int \Big(4r^{4N}\frac{\rho_{1,G}^4}{\rho_{2,G}^2}\Big[(2Z_{c_{61,i},1}-Z_{c_{61,i},2})^2c_{61,i}\Big]
-r^{4N}\frac{\rho_{2,G}}{\rho_{1,G}}\Big[(Z_{c_{61,i},2}-Z_{c_{61,i},1})(2Z_{c_{61,i},1}-Z_{c_{61,i},2})c_{61,i}
\Big]\\
&-&r^{2N}\frac{\rho_{1,G}^2}{\rho_{2,G}}(2\psi_{1,0}-\psi_{2,0})
\Big[(2Z_{c_{61,i},1}-Z_{c_{61,i},2})^2c_{61,i}\Big]\\
&+&4r^{4N}\frac{\rho_{2,G}^4}{\rho_{1,G}^6}\Big[\frac{1}{3}
(2Z_{c_{61,i},2}-3Z_{c_{61,i},1})^2c_{61,i}\Big]
-r^{4N}\frac{\rho_{2,G}}{\rho_{1,G}}\Big[(Z_{c_{61,i},2}-Z_{c_{61,i},1})(2Z_{c_{61,i},2}
-3Z_{c_{61,i},1})c_{61,i}
\Big]\\
&-&r^{2N}\frac{\rho_{2,G}^2}{\rho_{1,G}^3}(2\psi_{2,0}-3\psi_{1,0})\Big[\frac{1}{3}
(2Z_{c_{61,i},2}-3Z_{c_{61,i},1})^2c_{61,i}\Big]
\Big)rdrd\theta\\
&=&\int \Big[\Big(4r^{4N}\frac{\rho_{1,G}^4}{\rho_{2,G}^2}(2Z_{c_{61,i},1}-Z_{c_{61,i},2})^2
-r^{4N}\frac{\rho_{2,G}}{\rho_{1,G}}(Z_{c_{61,i},2}-Z_{c_{61,i},1})(2Z_{c_{61,i},1}-Z_{c_{61,i},2})\\
&+&\frac{4}{3}r^{4N}\frac{\rho_{2,G}^4}{\rho_{1,G}^6}(2Z_{c_{61,i},2}-3Z_{c_{61,i},1})^2
-r^{4N}\frac{\rho_{2,G}}{\rho_{1,G}}(Z_{c_{61,i},2}-Z_{c_{61,i},1})(2Z_{c_{61,i},2}-3Z_{c_{61,i},1})\Big)\\
&-&\frac{1}{3}r^{2N}\Big(\frac{\rho_{1,G}^2}{\rho_{2,G}}(2Z_{c_{61,i},1}-Z_{c_{61,i},2})^2
+\frac{1}{3}\frac{\rho_{2,G}^2}{\rho_{1,G}^3}(2Z_{c_{61,i},2}-3Z_{c_{61,i},1})^2\Big)\psi
\Big)c_{61,i}\\
&=&\pi(J_{7}+\int q_5 \psi r dr)c_{61,i}+O(\ve)+O(|\mathbf{a}|^2+\ve ^2),
\end{eqnarray*}
\begin{eqnarray*}
&&\int E\cdot Z_{c_{62,i}}^* rdrd\theta\\
&=&\int \Big(4r^{4N}\frac{\rho_{1,G}^4}{\rho_{2,G}^2}\Big[(2Z_{c_{62,i},1}-Z_{c_{62,i},2})^2c_{62,i}\Big]
-r^{4N}\frac{\rho_{2,G}}{\rho_{1,G}}\Big[(Z_{c_{62,i},2}-Z_{c_{62,i},1})(2Z_{c_{62,i},1}-Z_{c_{62,i},2})c_{62,i}
\Big]\\
&-&r^{2N}\frac{\rho_{1,G}^2}{\rho_{2,G}}(2\psi_{1,0}-\psi_{2,0})\Big[(2Z_{c_{62,i},1}-Z_{c_{62,i},2})^2
c_{62,i}\Big]\\
&+&4r^{4N}\frac{\rho_{2,G}^4}{\rho_{1,G}^6}\Big[\frac{1}{3}
(2Z_{c_{62,i},2}-3Z_{c_{62,i},1})^2c_{62,i}\Big]-r^{4N}\frac{\rho_{2,G}}{\rho_{1,G}}\Big[(Z_{c_{62,i},2}-Z_{c_{62,i},1})
(2Z_{c_{62,i},2}-3Z_{c_{62,i},1})c_{62,i}
\Big]\\
&-&r^{2N}\frac{\rho_{2,G}^2}{\rho_{1,G}^3}(2\psi_{2,0}-3\psi_{1,0})\Big[\frac{1}{3}
(2Z_{c_{62,i},2}-3Z_{c_{62,i},1})^2c_{62,i}\Big]
\Big)rdrd\theta\\
&=&\int \Big[\Big(4r^{4N}\frac{\rho_{1,G}^4}{\rho_{2,G}^2}(2Z_{c_{62,i},1}-Z_{c_{62,i},2})^2
-r^{4N}\frac{\rho_{2,G}}{\rho_{1,G}}(Z_{c_{62,i},2}-Z_{c_{62,i},1})(2Z_{c_{62,i},1}-Z_{c_{62,i},2})\\
&+&\frac{4}{3}r^{4N}\frac{\rho_{2,G}^4}{\rho_{1,G}^6}(2Z_{c_{62,i},2}-3Z_{c_{62,i},1})^2
-r^{4N}\frac{\rho_{2,G}}{\rho_{1,G}}(Z_{c_{62,i},2}-Z_{c_{62,i},1})(2Z_{c_{62,i},2}-3Z_{c_{62,i},1})\Big)\\
&-&\frac{1}{3}r^{2N}\Big(\frac{\rho_{1,G}^2}{\rho_{2,G}}(2Z_{c_{62,i},1}-Z_{c_{62,i},2})^2
+\frac{1}{3}\frac{\rho_{2,G}^2}{\rho_{1,G}^3}(2Z_{c_{62,i},2}-3Z_{c_{62,i},1})^2\Big)\psi
\Big)c_{62,i}\\
&=&\pi(J_{8}+\int q_6 \psi r dr)c_{62,i}+O(\ve)+O(|\mathbf{a}|^2+\ve ^2),
\end{eqnarray*}
where all the terms $O(\ve)\leq C\ve$, $O(\ve^2)\leq C\ve^2$, $O(|{\bf a}|^2)\leq C|{\bf a}|^2$ for some positive constant $C$ independent of ${\bf a}$ and $\ve$ provided that they are small enough.

So we get that
\begin{eqnarray*}
&&(\langle E, Z_{c_{43,1}}^*\rangle,\langle E, Z_{c_{43,2}}^*\rangle,\langle E, Z_{c_{52,1}}^*\rangle,\langle E, Z_{c_{52,2}}^*\rangle,\\
&& \, \, \, \, \, \, \, \, \, \, \langle E, Z_{c_{53,1}}^*\rangle,\langle E, Z_{c_{53,2}}^*\rangle,\langle E, Z_{c_{54,1}}^*\rangle,\langle E, Z_{c_{54,2}}^*\rangle,\\
&& \, \, \, \, \, \, \, \, \, \,\langle E, Z_{c_{61,1}}^*\rangle,\langle E, Z_{c_{61,2}}^*\rangle,
\langle E, Z_{c_{62,1}}^*\rangle,\langle E, Z_{c_{62,2}}^*\rangle)^t\\
&&=\tilde{\mathcal{T}}(\mathbf{a})+O(|\mathbf{a}|^2)+O(\ve),\nonumber
\end{eqnarray*}
and
\begin{equation}\label{Tb}
\tilde{\mathcal{T}}=\left(\begin{array}{cccccccccccc}
T_1&0&0&0&0&0&T_2&0&0&0&0&0\\
0&T_1&0&0&0&0&0&T_2&0&0&0&0\\
0&0&T_3&0&0&0&0&0&0&0&0&0\\
0&0&0&T_3&0&0&0&0&0&0&0&0\\
0&0&0&0&T_4&0&0&0&0&0&0&0\\
0&0&0&0&0&T_4&0&0&0&0&0&0\\
T_5&0&0&0&0&0&T_6&0&0&0&0&0\\
0&T_5&0&0&0&0&0&T_6&0&0&0&0\\
0&0&0&0&0&0&0&0&T_7&0&0&0\\
0&0&0&0&0&0&0&0&0&T_7&0&0\\
0&0&0&0&0&0&0&0&0&0&T_8&0\\
0&0&0&0&0&0&0&0&0&0&0&T_8\\
\end{array}
\right),
\end{equation}
where
\begin{eqnarray*}
&&T_1=(J_1+\int q_1\psi r dr),\ T_2=(J_2+\int q_7\psi r dr),\\
&&T_3=(J_5+\int q_2 \psi r dr), \ \ \ T_4=(J_6+\int q_3 \psi r dr),\\
&&T_5=(J_3+\int q_4 \psi r dr),\ \ \ T_6=(J_4+\int q_7 \psi r dr),\\
&&T_7=(J_{7}+\int q_5 \psi r dr),\ \ \ T_8=(J_{8}+\int q_6 \psi r dr)
\end{eqnarray*}

The determinant of the matrix $\tilde{\mathcal{T}}$ is
\begin{eqnarray*}
&&(T_1T_6-T_2T_5)^2T_3^3T_4^2T_7^2T_8^2\\
&&=(J_5+\int q_2 \psi r dr)^2(J_6+\int q_3 \psi r dr)^2(J_{7}+\int q_5 \psi r dr)^2
(J_{8}+\int q_6 \psi r dr)^2\\
&&\times [(J_1+\int q_1\psi r dr)(J_4+\int q_7 \psi r dr)-
(J_2+\int q_7\psi r dr)(J_3+\int q_4 \psi r dr)]^2
\end{eqnarray*}

Next we prove that the matrix $\tilde{\mathcal{T}}$ is non-degenerate, i.e. the determinant of $\tilde{\mathcal{T}} $ is nonzero. For this purpose, we need to calculate the integrals $J_1$ to $J_{8}$, and $\int q_1\psi r dr$ to $\int q_7\psi r dr$ But in the integrals, there is the function $\psi$ for which the expression is unknown.
In order to get rid of $\psi$, we use integration by parts. The key observation is that for any $\phi$ satisfying $\phi(\infty)=0$, we have
\begin{equation}\label{cancelpsib}
\int_0^\infty [(\Delta+\frac{8(N+1)^2r^{2N}}{(1+r^{2N+2})^2})\psi]\phi rdr
=\int_0^\infty[(\Delta+\frac{8(N+1)^2r^{2N}}{(1+r^{2N+2})^2})\phi]\psi rdr.
\end{equation}
By direct calculation, one can get that
\begin{eqnarray*}
q_1&&=-\frac{1}{3}r^{2N}[e^{2\tilde{U}_{1,0}-\tilde{U}_{2,0}}(2Z_{c_{43},1}(r)-Z_{c_{43},2}(r))^2
+\frac{1}{3}e^{2\tilde{U}_{2,0}-3\tilde{U}_{1,0}}(2Z_{c_{43},2}(r)-3Z_{c_{43},1}(r))^2]\\
&&=-18\frac{ \mu ^2 r^{2 \mu+2N }}{(r^{2 \mu }+1)^{12}}
(-178 r^{2 \mu }+1252 r^{4 \mu }-3746 r^{6 \mu }+5380 r^{8 \mu }-3746 r^{10 \mu }+1252 r^{12 \mu}\\
&&-178 r^{14 \mu }+9 r^{16 \mu }+9),\\
q_2&&=-\frac{1}{3}r^{2N}[e^{2\tilde{U}_{1,0}-\tilde{U}_{2,0}}(2Z_{c_{52},1}(r)-Z_{c_{52},2}(r))^2
+\frac{1}{3}e^{2\tilde{U}_{2,0}-3\tilde{U}_{1,0}}(2Z_{c_{52},2}(r)-3Z_{c_{52},1}(r))^2]\\
&&=-2240\frac{ \mu ^2 r^{6 \mu+2N }}{(r^{2 \mu }+1)^{12}}{ (-5 r^{2 \mu }+2 r^{4 \mu }+2)^2},\\
q_3&&=-\frac{1}{3}r^{2N}[e^{2\tilde{U}_{1,0}-\tilde{U}_{2,0}}(2Z_{c_{53},1}(r)-Z_{c_{53},2}(r))^2
+\frac{1}{3}e^{2\tilde{U}_{2,0}-3\tilde{U}_{1,0}}(2Z_{c_{53},2}(r)-3Z_{c_{53},1}(r))^2]\\
&&=-1260\frac{ \mu ^2 r^{4 \mu+2N } }{(r^{2 \mu }+1)^{12}}{(5 r^{2 \mu }-5 r^{4 \mu }+r^{6 \mu }-1)^2},\\
q_4&&=-\frac{1}{3}r^{2N}[e^{2\tilde{U}_{1,0}-\tilde{U}_{2,0}}(2Z_{c_{54},1}(r)-Z_{c_{54},2}(r))^2
+\frac{1}{3}e^{2\tilde{U}_{2,0}-3\tilde{U}_{1,0}}(2Z_{c_{54},2}(r)-3Z_{c_{54},1}(r))^2]\\
&&=-32\frac{ \mu ^2 r^{2 \mu+2N }}{9 (r^{2 \mu }+1)^{12}} \Big(-578 r^{2 \mu }+4052 r^{4 \mu }\\
&&-12146 r^{6 \mu }+17420 r^{8 \mu }-12146 r^{10 \mu }+4052 r^{12 \mu }-578 r^{14 \mu }+29 r^{16 \mu }+29\Big),
\\
q_5&&=-\frac{1}{3}r^{2N}[e^{2\tilde{U}_{1,0}-\tilde{U}_{2,0}}(2Z_{c_{61},1}(r)-Z_{c_{61},2}(r))^2
+\frac{1}{3}e^{2\tilde{U}_{2,0}-3\tilde{U}_{1,0}}(2Z_{c_{61},2}(r)-3Z_{c_{61},1}(r))^2]\\
&&=\frac{-560 \mu ^2 r^{10 \mu+2N }}{(r^{2 \mu }+1)^{12}},\\
q_6&&=-\frac{1}{3}r^{2N}[e^{2\tilde{U}_{1,0}-\tilde{U}_{2,0}}(2Z_{c_{62},1}(r)-Z_{c_{62},2}(r))^2
+\frac{1}{3}e^{2\tilde{U}_{2,0}-3\tilde{U}_{1,0}}(2Z_{c_{62},2}(r)-3Z_{c_{62},1}(r))^2]\\
&&=\frac{-560 \mu ^2 r^{8 \mu+2N } (r^{2 \mu }-1)^2}{(r^{2 \mu }+1)^{12}},\\
q_7&&=-\frac{1}{3}r^{2N}[e^{2\tilde{U}_{1,0}-\tilde{U}_{2,0}}(2Z_{c_{54},1}(r)-Z_{c_{54},2}(r))
(2Z_{c_{43},1}(r)-Z_{c_{43},2}(r))\\
&&
+\frac{1}{3}e^{2\tilde{U}_{2,0}-3\tilde{U}_{1,0}}(2Z_{c_{54},2}(r)-3Z_{c_{54},1}(r))(2Z_{c_{43},2}(r)
-3Z_{c_{43},1})(r)]\\
&&=16\frac{ \mu ^2 r^{2 \mu+2N } }{(r^{2 \mu }+1)^{12}}{(-11 r^{2 \mu }+21 r^{4 \mu }-11 r^{6 \mu }+r^{8 \mu }+1) (-73 r^{2 \mu }+153 r^{4 \mu }-73 r^{6 \mu }+8 r^{8 \mu }+8)},
\end{eqnarray*}
where we denote by $Z_{c_{ij},k}(r)$ the radial part of $Z_{c_{ij},k}$.

\medskip

So we need to find solutions of ODEs:
\begin{equation}
\Delta \phi_i+\frac{8(N+1)^2r^{2N}}{(1+r^{2N+2})^2}\phi_i=q_i,
\end{equation}
for $q_i$ defined as above and $i=1,\cdots,7$.

Following the same idea as in the proof of Lemma 3.4 in \cite{ALW}, one can get that
\begin{eqnarray*}
&&\phi_1=-\frac{ (81 r^{16\mu}-352 r^{14\mu}+1393 r^{12\mu}-1584 r^{10\mu}+1435 r^{8\mu}-304 r^{6\mu}+108 r^{4\mu}+8 r^{2\mu}+1)}{8 (r^{2\mu}+1)^{10}},\\
&&\phi_2=-\frac{2 (350 r^{12\mu}-336 r^{10\mu}+392 r^{8\mu}+48 r^{6\mu}+27 r^{4\mu}+8 r^{2\mu}+1)}{5 (r^{2\mu}+1)^{10}},\\
&&\phi_3=-\frac{ (350 r^{14\mu}-875 r^{12\mu}+1764 r^{10\mu}-833 r^{8\mu}+398 r^{6\mu}+27 r^{4\mu}+8 r^{2\mu}+1)}{10 (r^{2\mu}+1)^{10}},\\
&&\phi_4=-\frac{2 (145 r^{16\mu}-640 r^{14\mu}+2485 r^{12\mu}-2896 r^{10\mu}+2527 r^{8\mu}-592 r^{6\mu}+172 r^{4\mu}+8 r^{2\mu}+1)}{45 (r^{2\mu}+1)^{10}},\\
&&\phi_5=-\frac{5 (42 r^{8\mu}+48 r^{6\mu}+27 r^{4\mu}+8 r^{2\mu}+1)}{54 (r^{2\mu}+1)^{10}},\\
&&\phi_6=-\frac{4 (189 r^{10\mu}+42 r^{8\mu}+48 r^{6\mu}+27 r^{4\mu}+8 r^{2\mu}+1)}{135 (r^{2\mu}+1)^{10}},\\
&&\phi_7=\frac{ r^{4\mu} (16 r^{12\mu}-72 r^{10\mu}+273 r^{8\mu}-328 r^{6\mu}+273 r^{4\mu}-72 r^{2\mu}+16)}{2 (r^{2\mu}+1)^{10}}.
\end{eqnarray*}

So by direct calculation, we have
\begin{eqnarray*}
\int\psi q_i rdr=\int_0^\infty r^{4N}\frac{3\times 2^6 (N+1)^4}{(1+r^{2N+2})^4}\phi_i rdr.
\end{eqnarray*}

Since all the terms in the integrals are explicit now, by direct calculation, we get

\begin{eqnarray*}
&&J_5+\int q_2\psi r dr\\
&&=\frac{\pi  N\csc(\frac{\pi}{N+1})}{5405400 (N+1)^{10}}\\
&&\ \ \ \ \times \Big(229219200 N^{11}+2741981760 N^{10}+14845935288 N^9+48087228720 N^8\\
&&\ \ \ \ \ \ \ +103607697806 N^7+155921208688 N^6+167142190971 N^5+127469747650 N^4\\
&&\ \ \ \ \ \ \ \ +67655345084 N^3 +23742751992 N^2+4943660256 N+461099520\Big),\\
&&J_6+\int q_3\psi r dr\\
&&=\frac{N \pi \csc(\frac{\pi}{N+1})}{5405400 (N+1)^{10}}\\
&& \ \ \ \ \times \Big(73483200 N^{11}+879636240 N^{10}+4780746072 N^9+15595933680 N^8+33964635664 N^7\\
&&\ \ \ \ \ \ \ +51871736672 N^6+56687585199 N^5+44305756250 N^4+24240925096 N^3\\
&&\ \ \ \ \ \ \ \ +8824079448 N^2+1917285264 N+187548480\Big),\\
\end{eqnarray*}
\begin{eqnarray*}
&&J_7+\int q_5\psi r dr\\
&&=\frac{\pi  N\csc(\frac{\pi}{N+1}) }{70053984 (N+1)^{10}}\\
&&\ \ \ \ \times
\Big(604195200 N^{11}+7224396480 N^{10}+38670644088 N^9+122313524400 N^8\\
&&\ \ \ \ \ \ \ +253958797454 N^7+363332551792 N^6+365322970803 N^5+257996555026 N^4\\
&&\ \ \ \ \ \ \ \ +125308378700 N^3+39819491064 N^2+7439155488 N+617621760\Big),
\\
&&J_8+\int q_6\psi r dr\\
&&=\frac{\pi  N \csc (\frac{\pi }{N+1})}{437837400 (N+1)^{10}}\\
&&\ \ \ \ \times \Big((N+2) (2 N+3) (3 N+4) (4 N+5) (5 N+6) (N (N (2 N (12 N (504 N (855 N+3997)\\
&&\ \ \ \ \ \ \ +4024843)+51916217)+62504971)+19787638)+2554776)\Big)
,
\end{eqnarray*}
\begin{eqnarray*}
&&(J_1+\int q_1\psi r dr )(J_4+\int q_7\psi r dr)-(J_2+\int q_7\psi r dr )(J_3+\int q_4\psi r dr)\\
&&=\frac{\pi ^2 N^3  \csc ^2(\frac{\pi }{N+1})}{6949800 (N+1)^{12}}
\Big((N+2)^2 (2 N+3)^2 (N (N (N (N (2 N (9 N (2 N (120 N (4365 N+37031)\\
&&\ \ \ \ \ \ \ +17041651)+77384503)+1031119715)+2064834729)+1398543708)\\
&&\ \ \ \ \ \ \ \ +617719460)+161120544)+18780480)\Big).
\end{eqnarray*}

One can check that the above terms are all nonzero if $N>0$. So $\tilde{\mathcal{T}}$ is non-degenerate if $N>0$.

\qed

As a consequence of Lemma \ref{lemma401} and Lemma \ref{lemma402b}, the coefficients $m_i$ varnish if and only if the parameters ${\bf a}$
satisfy
\begin{eqnarray}\label{a}
\tilde{\mathcal{T}}(\mathbf{a})+O(|\mathbf{a}|^2)+O(\ve)=0,
\end{eqnarray}
where $\tilde{\mathcal{T}}$ is defined in (\ref{Tb}). Obviously, (\ref{a}) can be solved immediately from (\ref{Tb}) with $|\mathbf{a}|\leq C\ve$, for some $C$ large but fixed.

\medskip

If $N=0$, we choose the zero-th approximate solutions to be radial $\tilde{U}_0$, and the final approximate solution to be radial, then one can easily find a radial solution $\vect{v_1}{v_2}$ to (\ref{projectb}) with $m_i$ all zero in Proposition \ref{propb}, since $G$ is radial, (\ref{errorb}) is automatically satisfied.

\section{Proof of Theorem {\ref{thm101}} under Assumption (ii)}\label{sec3.6}

In this section, we are going to prove Theorem \ref{thm101} under {\bf Assumption (ii)}. This situation is more complicated than the previous one, since the $O(\ve)$ approximation and $O(\ve^2)$ approximation induce several difficulties.  The problem is that we cannot obtain the explicit expressions for these terms. In this case, we will see that the two free parameters $\xi_1$, $\xi_2$ we introduced in Section \ref{sec3.4} for the improvement of the $O(\ve^2)$ approximate solution play an important role. A key observation is that we only need to consider the terms involving $\xi_1$ and $\xi_2$. This is contained in the following lemma.

\begin{lemma}\label{lemma501b}
Let $\vect{v_1}{v_2}$ be a solution of (\ref{projectb}). The following estimates hold:
\begin{eqnarray}
&&\langle E, Z_{c_{43,i}}^*\rangle \\
&&=\xi_1(\tilde{\mathcal{A}}_1c_{43,i}+\tilde{\mathcal{B}}_1c_{54,i})
+\xi_2(\tilde{\mathcal{A}}_2c_{43,i}+
\tilde{\mathcal{B}}_2c_{54,i})+\tilde{\mathcal{T}}_{1i}(\mathbf{a})
+O((1+|\xi|)|\mathbf{a}|^2)+O(\ve),\nonumber\\
&&\langle E, Z_{c_{52,i}}^*\rangle \\
&&=\xi_1\tilde{\mathcal{C}}_1c_{52,i}+\xi_2\tilde{\mathcal{C}}_2c_{52,i}
+\tilde{\mathcal{T}}_{3i}(\mathbf{a})+O((1+|\xi|)|\mathbf{a}|^2)
+O(\ve),\nonumber\\
&&\langle E, Z_{c_{53,i}}^*\rangle \\
&&=\xi_1\tilde{\mathcal{D}}_1c_{53,i}+\xi_2\tilde{\mathcal{D}}_2c_{53,i}
+\tilde{\mathcal{T}}_{4i}(\mathbf{a})+O((1+|\xi|)|\mathbf{a}|^2)
+O(\ve),\nonumber\\
&&\langle E, Z_{c_{54,i}}^*\rangle \\
&&=\xi_1(\tilde{\mathcal{E}}_1c_{54,i}+\tilde{\mathcal{F}}_1c_{43,i})
+\xi_2(\tilde{\mathcal{E}}_2c_{54,i}+
\tilde{\mathcal{F}}_2c_{43,i})+\tilde{\mathcal{T}}_{2i}(\mathbf{a})
+O((1+|\xi|)|\mathbf{a}|^2)+O(\ve),\nonumber\\
&&\langle E, Z_{c_{61,i}}^*\rangle \\
&&=\xi_1\tilde{\mathcal{G}}_1c_{61,i}+\xi_2\tilde{\mathcal{G}}_2c_{61,i}
+\tilde{\mathcal{T}}_{4i}(\mathbf{a})+O((1+|\xi|)|\mathbf{a}|^2)
+O(\ve),\nonumber\\
&&\langle E, Z_{c_{62,i}}^*\rangle \\
&&=\xi_1\tilde{\mathcal{H}}_1c_{62,i}+\xi_2\tilde{\mathcal{H}}_2c_{62,i}
+\tilde{\mathcal{T}}_{4i}(\mathbf{a})+O((1+|\xi|)|\mathbf{a}|^2)
+O(\ve),\nonumber
\end{eqnarray}
for $i=1,2$, where
\begin{eqnarray}\label{a1b}
\tilde{\mathcal{A}}_j&=&\int_0^{+\infty}\int_0^{2\pi} \Big[ r^{2N_1}e^{2\tilde{U}_{1,0}-\tilde{U}_{2,0}}(2Z_{j,1}-Z_{j,2})(2Z_{c_{43,1},1}-Z_{c_{43,1},2})^2\nonumber\\
&+&\frac{1}{3}r^{2N_2}e^{2\tilde{U}_{2,0}-3\tilde{U}_{1,0}}(2Z_{j,2}-3Z_{j,1})
(2Z_{c_{43,1},2}-3Z_{c_{43,1},1})^2 \Big]r drd\theta,
\end{eqnarray}
\begin{eqnarray}\label{b1b}
&&\tilde{\mathcal{B}}_j=\int_0^{+\infty}\int_0^{2\pi} \Big[ r^{2N_1}e^{2\tilde{U}_{1,0}-\tilde{U}_{2,0}}(2Z_{j,1}-Z_{j,2})(2Z_{c_{43,1},1}-Z_{c_{43,1},2})
(2Z_{c_{54,1},1}-Z_{c_{54,1},2})\nonumber\\
&&+\frac{1}{3}r^{2N_2}e^{2\tilde{U}_{2,0}-3\tilde{U}_{1,0}}(2Z_{j,2}-3Z_{j,1})
(2Z_{c_{54,1},2}-3Z_{c_{54,1},1})(2
Z_{c_{43,1},2}-3Z_{c_{43,1},1})\Big]r drd\theta,\nonumber\\
\end{eqnarray}

\begin{eqnarray}
\tilde{\mathcal{C}}_j&=&\int_0^{+\infty}\int_0^{2\pi} \Big[ r^{2N_1}e^{2\tilde{U}_{1,0}-\tilde{U}_{2,0}}(2Z_{j,1}-Z_{j,2})
(2Z_{c_{52,1},1}-Z_{c_{52,1},2})^2\nonumber\\
&+&\frac{1}{3}r^{2N_2}e^{2\tilde{U}_{2,0}-3\tilde{U}_{1,0}}
(2Z_{j,2}-3Z_{j,1})(2Z_{c_{52,1},2}-3Z_{c_{52,1},1})^2\Big] r drd\theta,
\end{eqnarray}

\begin{eqnarray}
\tilde{\mathcal{D}}_j&=&\int_0^{+\infty}\int_0^{2\pi} \Big[ r^{2N_1}e^{2\tilde{U}_{1,0}-\tilde{U}_{2,0}}(2Z_{j,1}-Z_{j,2})
(2Z_{c_{53,1},1}-Z_{c_{53,1},2})^2\nonumber\\
&+&\frac{1}{3}r^{2N_2}e^{2\tilde{U}_{2,0}-3\tilde{U}_{1,0}}
(2Z_{j,2}-3Z_{j,1})(2Z_{c_{53,1},2}-3Z_{c_{53,1},1})^2\Big] r drd\theta,
\end{eqnarray}

\begin{eqnarray}
\tilde{\mathcal{E}}_j&=&\int_0^{+\infty}\int_0^{2\pi} \Big[ r^{2N_1}e^{2\tilde{U}_{1,0}-\tilde{U}_{2,0}}(2Z_{j,1}-Z_{j,2})
(2Z_{c_{54,1},1}-Z_{c_{54,1},2})^2\nonumber\\
&+&\frac{1}{3}r^{2N_2}e^{2\tilde{U}_{2,0}-3\tilde{U}_{1,0}}(2Z_{j,2}-3Z_{j,1})
(2Z_{c_{54,1},2}-3Z_{c_{54,1},1})^2\Big] r drd\theta,
\end{eqnarray}
\begin{eqnarray}
&&\tilde{\mathcal{F}}_j=\int_0^{+\infty}\int_0^{2\pi} \Big[ r^{2N_1}e^{2\tilde{U}_{1,0}-\tilde{U}_{2,0}}(2Z_{j,1}-Z_{j,2})(2Z_{c_{43,1},1}-Z_{c_{43,1},2})
(2Z_{c_{54,1},1}-Z_{c_{54,1},2})\nonumber\\
&&+\frac{1}{3}r^{2N_2}e^{2\tilde{U}_{2,0}-3\tilde{U}_{1,0}}(2Z_{j,2}-3Z_{j,1})
(2Z_{c_{54,1},2}-3Z_{c_{54,1},1})(2
Z_{c_{43,1},2}-3Z_{c_{43,1},1})\Big]r drd\theta,\nonumber\\
\end{eqnarray}

\begin{eqnarray}
\tilde{\mathcal{G}}_j&=&\int_0^{+\infty}\int_0^{2\pi} \Big[ r^{2N_1}e^{2\tilde{U}_{1,0}-\tilde{U}_{2,0}}(2Z_{j,1}-Z_{j,2})(2Z_{c_{61,1},1}-Z_{c_{61,1},2})^2\nonumber\\
&+&\frac{1}{3}r^{2N_2}e^{2\tilde{U}_{2,0}-3\tilde{U}_{1,0}}(2Z_{j,2}-3Z_{j,1})
(2Z_{c_{61,1},2}-3Z_{c_{61,1},1})^2\Big] r drd\theta,
\end{eqnarray}
\begin{eqnarray}
\tilde{\mathcal{H}}_j&=&\int_0^{+\infty}\int_0^{2\pi} \Big[ r^{2N_1}e^{2\tilde{U}_{1,0}-\tilde{U}_{2,0}}(2Z_{j,1}-Z_{j,2})
(2Z_{c_{62,1},1}-Z_{c_{62,1},2})^2\nonumber\\
&+&\frac{1}{3}r^{2N_2}e^{2\tilde{U}_{2,0}-3\tilde{U}_{1,0}}(2Z_{j,2}-3Z_{j,1})
(2Z_{c_{62,1},2}-3Z_{c_{62,1},1})^2 \Big]r drd\theta,
\end{eqnarray}
for $j=1,2$, and $\tilde{\mathcal{T}}_{ij}$ are $12\times 1$ vectors which are uniformly bounded as $\ve$ tends to $0$, and are independent of $\xi_1,\xi_2$.
\end{lemma}
\noindent
{\bf Proof:}

By \eqref{E1000b} and \eqref{E2000b}, $E$ is of the form
\begin{equation}
\label{E100b}
\frac{1}{\ve}   (...) {\bf a} \cdot {\bf a} +  ( (...) {\bf a} ) + O(  |{\bf a}|^2) +O(\ve).
\end{equation}

Recall that $\vect{\psi_1}{\psi_2}=\vect{\psi_{0,1}}{\psi_{0,2}}
+\xi_1Z_{\lambda_4}+\xi_2Z_{\lambda_5}$.
In the following computations, we only need to consider  the terms involving $\xi_1$ and $\xi_2$, since all other terms are independent of $\xi_1$ and $\xi_2$.  By the orthogonality of $\cos(k\theta)$ and $\cos(l\theta)$ for $k\neq l$, we obtain
\begin{eqnarray*}
&-&\int_0^{+\infty}\int_0^{2\pi} E\cdot Z_{c_{43,1}}^*rd\theta dr\\
&=&\int_0^\infty\int_0^{2\pi} \Big[r^{2N_1}e^{2\tilde{U}_{1,0}-\tilde{U}_{2,0}}\Big(\xi_1(2Z_{1,1}-Z_{1,2})
+\xi_2(2Z_{2,1}-Z_{2,2})\Big)\\
&\times &(2Z_{c_{43,1},1}-Z_{c_{43,1},2})c_{43,1}Z_{c_{43,1},1}^*\\
&+&r^{2N_2}e^{2\tilde{U}_{2,0}-3\tilde{U}_{1,0}}\Big(\xi_1(2Z_{1,2}-3Z_{1,1})
+\xi_2(2Z_{2,2}-3Z_{2,1})\Big)\\
&\times &(2Z_{c_{43,1},2}-3Z_{c_{43,1},1})c_{43,1}Z_{c_{43,1},2}^*\\
&+&r^{2N_1}e^{2\tilde{U}_{1,0}-\tilde{U}_{2,0}}\Big(\xi_1(2Z_{1,1}-Z_{1,2})
+\xi_2(2Z_{2,1}-Z_{2,2})\Big)\\
&\times &(2Z_{c_{54,1},1}-Z_{c_{54,1},2})c_{54,1}Z_{c_{43,1},1}^*\\
&+&r^{2N_2}e^{2\tilde{U}_{2,0}-3\tilde{U}_{1,0}}\Big(\xi_1(2Z_{1,2}-3Z_{1,1})
+\xi_2(2Z_{2,2}-3Z_{2,1})\Big)\\
&\times &(2Z_{c_{54,1},2}-3Z_{c_{54,1},1})c_{54,1}Z_{c_{43,1},2}^*\Big]rdrd\theta\\
&+&\tilde{\mathcal{T}}_{11}(\mathbf{a})+O(\frac{|{\bf a}|^2}{\ve})+O(\ve)+O((1+|\xi|)|{\bf a}|^2+\ve ^2),
\end{eqnarray*}
where $\tilde{\mathcal{T}}_{11}(\mathbf{a})$ is the remaining terms which is a  linear combinations of $\mathbf{a}$ which comes from the remaining terms of $O(\ve^2)$ of $E$, {and the coefficients of the linear combinations are uniformly bounded and are independent of $\xi_1,\xi_2,\mathbf{a} $. The $O(\frac{|{\bf a}|^2}{\ve})$ terms comes from the $O(\frac{|{\bf a}|^2}{\ve})$ term of $E$ which is independent of $\xi$}.

Thus
\begin{eqnarray*}
&&-\int _0^{+\infty }\int _0^{2\pi }E\cdot Z_{c_{43,1}}^*rd\theta dr \\
&&=\xi_1\Big[\Big(\int_0^{+\infty}\int_0^{2\pi} \Big\{ r^{2N_1}e^{2\tilde{U}_{1,0}-\tilde{U}_{2,0}}(2Z_{1,1}-Z_{1,2})(2Z_{c_{43,1},1}-Z_{c_{43,1},2})^2\\
&&+\frac{1}{3}r^{2N_2}e^{2\tilde{U}_{2,0}-3\tilde{U}_{1,0}}(2Z_{1,2}-3Z_{1,1})
(2Z_{c_{43,1},2}-3Z_{c_{43,1},1})^2\Big\} r drd\theta\Big)c_{43,1}\\
&&+\Big(\int_0^{+\infty}\int_0^{2\pi} \Big\{ r^{2N_1}e^{2\tilde{U}_{1,0}-\tilde{U}_{2,0}}(2Z_{1,1}-Z_{1,2})(2Z_{c_{54,1},1}-Z_{c_{54,1},2})
(2Z_{c_{43,1},1}-Z_{c_{43,1},2})\\
&&+\frac{1}{3}r^{2N_2}e^{2\tilde{U}_{2,0}-3\tilde{U}_{1,0}}(2Z_{1,2}-3Z_{1,1})
(2Z_{c_{43,1},2}-3Z_{c_{43,1},1})
(2Z_{c_{54,1},2}-3Z_{c_{54,1},1})\Big\}r drd\theta\Big)c_{54,1}\Big]\\
&&+\xi_2\Big[\Big(\int_0^{+\infty}\int_0^{2\pi} \Big\{ r^{2N_1}e^{2\tilde{U}_{1,0}-\tilde{U}_{2,0}}(2Z_{2,1}-Z_{2,2})(2Z_{c_{43,1},1}-Z_{c_{43,1},2})^2\\
&&+\frac{1}{3}r^{2N_2}e^{2\tilde{U}_{2,0}-3\tilde{U}_{1,0}}(2Z_{2,2}-3Z_{2,1})
(Z_{c_{43,1},2}-Z_{c_{43,1},1})^2\Big\} r drd\theta\Big)c_{43,1}\\
&&+\Big(\int_0^{+\infty}\int_0^{2\pi}\Big\{ r^{2N_1}e^{2\tilde{U}_{1,0}-\tilde{U}_{2,0}}(2Z_{2,1}-Z_{2,2})(2Z_{c_{43,1},1}-Z_{c_{43,1},2})
(2Z_{c_{54,1},1}-Z_{c_{54,1},2})\\
&&+\frac{1}{3}r^{2N_2}e^{2\tilde{U}_{2,0}-3\tilde{U}_{1,0}}
(2Z_{2,2}-3Z_{2,1})(2Z_{c_{43,1},2}-3Z_{c_{43,1},1})
(2Z_{c_{54,1},2}-3Z_{c_{54,1},1})\Big\}r drd\theta\Big)c_{54,1}\Big]\\
&&+\tilde{\mathcal{T}}_{11}(\mathbf{a})+O(\frac{|{\bf a}|^2}{\ve})+O(\ve)+O((1+|\xi|)|{\bf a}|^2+\ve ^2)\\
&&=\xi_1(\tilde{\mathcal{A}}_1c_{43,1}+\tilde{\mathcal{B}}_1c_{54,1})
+\xi_2(\tilde{\mathcal{A}}_2c_{43,1}+\tilde{\mathcal{B}}_2c_{54,1})\\
&&+\tilde{\mathcal{T}}_{11}(\mathbf{a})+O(\frac{|{\bf a}|^2}{\ve})+O(\ve)++O((1+|\xi|)|{\bf a}|^2+\ve ^2),
\end{eqnarray*}
where $\tilde{\mathcal{A}}_1,\tilde{\mathcal{A}}_2,\tilde{\mathcal{B}}_1,\tilde{\mathcal{B}}_2$ are in (\ref{a1b}) and (\ref{b1b}).

Similarly, we can get the other estimates.

\qed

From the above lemma, we have the following result:
\begin{lemma}\label{lemma502b}
Let $\vect{v_1}{v_2}$ be a solution of (\ref{projectb}). Then the coefficients $m_i=0$ if and only if the parameters ${\bf a}$
satisfy
\begin{eqnarray}\label{ab}
\tilde{{\bf Q}}(\mathbf{a})=\tilde{\mathcal{T}}(\mathbf{a})+O(\frac{|{\bf a}|^2}{\ve})+O((1+|\xi|)|\mathbf{a}|^2)+O(\ve),
\end{eqnarray}
where $\tilde{\bf Q}$
\begin{eqnarray}\label{Qb}
\tilde{{\bf Q}}&=&\xi_1\left(\begin{array}{cccccccccccc}
\tilde{\mathcal{A}}_1&0&0&0&0&0&\tilde{\mathcal{B}}_1&0&0&0&0&0\\
0&\tilde{\mathcal{A}}_1&0&0&0&0&0&\tilde{\mathcal{B}}_1&0&0&0&0\\
0&0&\tilde{\mathcal{C}}_1&0&0&0&0&0&0&0&0&0\\
0&0&0&\tilde{\mathcal{C}}_1&0&0&0&0&0&0&0&0\\
0&0&0&0&\tilde{\mathcal{D}}_1&0&0&0&0&0&0&0\\
0&0&0&0&0&\tilde{\mathcal{D}}_1&0&0&0&0&0&0\\
\tilde{\mathcal{F}}_1&0&0&0&0&0&\tilde{\mathcal{E}}_1&0&0&0&0&0\\
0&\tilde{\mathcal{F}}_1&0&0&0&0&0&\tilde{\mathcal{E}}_1&0&0&0&0\\
0&0&0&0&0&0&0&0&\tilde{\mathcal{G}}_1&0&0&0\\
0&0&0&0&0&0&0&0&0&\tilde{\mathcal{G}}_1&0&0\\
0&0&0&0&0&0&0&0&0&0&\tilde{\mathcal{H}}_1&0\\
0&0&0&0&0&0&0&0&0&0&0&\tilde{\mathcal{H}}_1
\end{array}
\right)\\
&+&\xi_2\left(\begin{array}{cccccccccccc}
\tilde{\mathcal{A}}_2&0&0&0&0&0&\tilde{\mathcal{B}}_2&0&0&0&0&0\\
0&\tilde{\mathcal{A}}_2&0&0&0&0&0&\tilde{\mathcal{B}}_2&0&0&0&0\\
0&0&\tilde{\mathcal{C}}_2&0&0&0&0&0&0&0&0&0\\
0&0&0&\tilde{\mathcal{C}}_2&0&0&0&0&0&0&0&0\\
0&0&0&0&\tilde{\mathcal{D}}_2&0&0&0&0&0&0&0\\
0&0&0&0&0&\tilde{\mathcal{D}}_2&0&0&0&0&0&0\\
\tilde{\mathcal{F}}_2&0&0&0&0&0&\tilde{\mathcal{E}}_2&0&0&0&0&0\\
0&\tilde{\mathcal{F}}_2&0&0&0&0&0&\tilde{\mathcal{E}}_2&0&0&0&0\\
0&0&0&0&0&0&0&0&\tilde{\mathcal{G}}_2&0&0&0\\
0&0&0&0&0&0&0&0&0&\tilde{\mathcal{G}}_2&0&0\\
0&0&0&0&0&0&0&0&0&0&\tilde{\mathcal{H}}_2&0\\
0&0&0&0&0&0&0&0&0&0&0&\tilde{\mathcal{H}}_2
\end{array}
\right) \nonumber\\
&=&\xi_1{\bf Q}_1+\xi_2{\bf Q}_2,
\end{eqnarray}
and $\tilde{\mathcal{T}}$ is a $12\times 12$ matrix which is uniformly bounded and independent of $\xi_1,\xi_2$.
\end{lemma}

\medskip
\noindent
{\bf
Proof of Theorem \ref{thm101} under Assumption (ii):}
 \medskip
Under the {\bf Assumptions (ii)}, we will choose $\lambda_4=\frac{\tilde{\lambda}}{(2^{7/2} \mu_1 \mu_2 (\mu_1+\mu_2))^2}$ and $\lambda_5=\frac{1}{\tilde{\lambda}}$ for $\tilde{\lambda}$ large enough. Similar to the computation in the appendix of \cite{ALW}, by direct but tedious computation,  we can get that for $\tilde{\lambda}$ large,
\begin{eqnarray*}
&&\tilde{\mathcal{A}}_1=\gamma_1 \tilde{\lambda}^{-1}+o(\tilde{\lambda}^{-1}),\\
&&\tilde{\mathcal{C}}_1=\gamma_2 \tilde{\lambda}^{-1}+o(\tilde{\lambda}^{-1}),\\
&&\tilde{\mathcal{D}}_1=\gamma_3 \tilde{\lambda}^{-2}+o(\tilde{\lambda}^{-2}),\\
&&\tilde{\mathcal{E}}_1=\gamma_4 \tilde{\lambda}^{-4}+o(\tilde{\lambda}^{-4}),\\
&&\tilde{\mathcal{G}}_1=\gamma_5 \tilde{\lambda}^{-3}+o(\tilde{\lambda}^{-3}),\\
&&\tilde{\mathcal{H}}_1=\gamma_6 +o(1),\\
\end{eqnarray*}
and
\begin{equation}
\tilde{\mathcal{B}}_1=\tilde{\mathcal{F}}_1=
\left\{\begin{array}{c}
O(\tilde{\lambda}^{-3}) \mbox{ if } N_1=N_2,\\
0 \mbox{ if } N_1\neq N_2,
\end{array}
\right.
\end{equation}
where $\gamma_1$ to $\gamma_6$ are non zero as follows:
\begin{eqnarray*}
\gamma_1&&=\frac{4092 \mu_1 (\mu_1+\mu_2) (3 \mu_1+\mu_2)}{2 \mu_1+\mu_2}\pi^2 \csc (\frac{\mu_1\pi}{2\mu_1+\mu_2}),\\
\gamma_2&&=-44695552 \mu_1^4 \mu_2^2 (\mu_1+\mu_2)^2 (2 \mu_1+\mu_2)^3 (3 \mu_1+\mu_2)^2\pi,\\
\gamma_3&&=\mu_1^3 \mu_2^2 2^{36-\frac{56 \mu_1}{2 \mu_1+\mu_2}} (2 \mu_1-\mu_2) (6 \mu_1-\mu_2) (\mu_1+\mu_2)^4 (2 \mu_1+\mu_2) (3 \mu_1+2 \mu_2)\\
 &&\times \Big((\mu_1+\mu_2)^2 (2 \mu_1+\mu_2)^2 (3 \mu_1+2 \mu_2)^2\Big)^{1-\frac{8 \mu_1}{2 \mu_1+\mu_2}} (2 \mu_1^2-5 \mu_1 \mu_2-6 \mu_2^2)\pi^2 \csc (\frac{8\mu_1\pi}{2\mu_1+\mu_2}),\\
\gamma_4&&=\frac{\mu_1^4 \mu_2^5 2^{\frac{14 \mu_1}{2 \mu_1+\mu_2}+17}}{(3 \mu_1+\mu_2)^2} (\mu_1+\mu_2)^6 (2 \mu_1+\mu_2) ((\mu_1+\mu_2)^2 (2 \mu_1+\mu_2)^2 (3 \mu_1+2 \mu_2)^2)^{-\frac{2 (\mu_1+\mu_2)}{2 \mu_1+\mu_2}}\\
&&\times (9 \mu_1+7 \mu_2)\pi^2 \csc (\frac{2\mu_1\pi}{2\mu_1+\mu_2}),\\
\gamma_5&&=-341  \mu_1^6 \mu_2^6 2^{\frac{7 \mu_1}{2 \mu_1+\mu_2}+23} (\mu_1+\mu_2)^9 (3 \mu_1+2 \mu_2)^3 \\
&&\times ((\mu_1+\mu_2)^2 (2 \mu_1+\mu_2)^2 (3 \mu_1+2 \mu_2)^2)^{\frac{\mu_1}{2 \mu_1+\mu_2}-1} \times (5 \mu_1+3 \mu_2) \pi^2\csc (\frac{\pi  \mu_1}{2 \mu_1+\mu_2}),\\
\gamma_6&&=-\frac{1}{3}   \mu_1^4 \mu_2^2 2^{\frac{7 \mu_1}{2 \mu_1+\mu_2}+15} (\mu_1+\mu_2)^5 (2 \mu_1+\mu_2) (3 \mu_1+\mu_2) (3 \mu_1+2 \mu_2)^2 \\
&&\times ((\mu_1+\mu_2) (2 \mu_1+\mu_2) (3 \mu_1+2 \mu_2)) ^{\frac{2 \mu_1}{2 \mu_1+\mu_2}-4} \pi^2\csc (\frac{\pi  \mu_1}{2 \mu_1+\mu_2}).
\end{eqnarray*}

It is easy to see that from the above expressions that $\gamma_1,\gamma_2, \gamma_4,\gamma_5,\gamma_6$ are all non zero, and $\gamma_3$ is also non zero, since for $\mu_2=2\mu_1$ or $\mu_2=6\mu_1$, $\gamma_3$ is also non zero.

 \noindent
So we have
\begin{eqnarray*}
\tilde{\mathcal{A}}_1\tilde{\mathcal{E}}_1-\tilde{\mathcal{B}}_1\tilde{\mathcal{F}}_1&=&
\gamma_1\gamma_4\tilde{\lambda}^{-5}+o(\tilde{\lambda}^{-5})\\
&\neq & 0,
\end{eqnarray*}
and $\tilde{\mathcal{C}}_1, \tilde{\mathcal{D}}_1, \tilde{\mathcal{E}}_1, \tilde{\mathcal{G}}_1, \tilde{\mathcal{H}}_1$ are both non-zero if $\tilde{\lambda}$ is large enough.
Therefore, we choose $\xi_1$ large and $\xi_2=0$ to conclude that ${\bf Q}(\xi_1,\xi_2)-\tilde{\mathcal{T}}$ is non-degenerate. After fixing $({\lambda}_4,\lambda_5)$, $(\xi_1,\xi_2)$, it is easy to see (\ref{a1b}) can be solved with $\mathbf{a}=O(\ve)$.

\qed

\section{Proof of Theorem \ref{thm101} under Assumption (iii) }\label{sec3.7}

We are left to prove the theorem for $N_2\sum_{i=1}^{N_1}p_i=N_1\sum_{j=1}^{N_2}q_j$, $N_1\neq N_2$ and one of $N_i$ is $1$. Without loss of generality, assume $N_1=1$ and $\sum_{i=1}^{N_1}p_i=\sum_{j=1}^{N_2}q_j=0$. In this case, for the improvement of approximate solution in the $O(\ve^2)$ term, we can not solve equation (\ref{psib}) in Section \ref{sec3.4}. Instead of solving (\ref{psib}),  we can find a unique solution of the following equations which is guaranteed by Lemma \ref{lemma2b}:

\begin{equation}\label{psi2b}
\left\{\begin{array}{c}
\Delta \psi_1+|z|^{2N_1}e^{2\tilde{U}_{1,0}-\tilde{U}_{2,0}}(2\psi_1-\psi_2)=
2|z|^{4N_1}e^{4\tilde{U}_{1,0}-2\tilde{U}_{2,0}}-|z|^{2(N_1+N_2)}e^{\tilde{U}_{2,0}-\tilde{U}_{1,0}}\\
\Delta \psi_2+|z|^{2N_2}e^{2\tilde{U}_{2,0}-3\tilde{U}_{1,0}}(2\psi_2-3\psi_1)=
2|z|^{4N_2}e^{4\tilde{U}_{2,0}-6\tilde{U}_{1,0}}-3|z|^{2(N_1+N_2)}e^{\tilde{U}_{2,0}-\tilde{U}_{1,0}}.
\end{array}
\right.
\end{equation}

We use this unique solution as the new $\psi_0$, and proceed as before. Then by checking the previous proof, we can get that in this case, the error $\|E\|_*\leq C_0$ and we can get a solution $v$ of (\ref{projectb}) which satisfies
\begin{equation}
\|v\|_*\leq C_0,
\end{equation}
for some positive constant $C_0$, and the following estimates hold:

\begin{equation}
\int_{\R^2} (N_{11}(v)+N_{12}(v))Z_{i,1}^*+(N_{21}(v)+N_{22}(v))Z_{i,2}^* dx=O(\ve),
\end{equation}
for $i=3,\cdots,14$.

Then the reduced problem we get is
\begin{equation}\label{specialb}
{{\bf Q}}(\mathbf{a})+\tilde{\mathcal{T}}(\mathbf{a})+O((1+|\xi|)|\mathbf{a}|^2)+O(1)+O(\ve)=0,
\end{equation}
where the $O(1)$ term comes from the $O(1)$ term of the error $E$ since we use the solution of (\ref{psi2b}) instead of (\ref{psib}) as the $O(\ve^2)$ improvement. Recalling that ${\bf Q}={\bf Q}(\xi_1,\xi_2)$ depend on two free parameters $\xi_1,\xi_2$ and arguing as before, we can choose $\xi_1$ large enough. Then it is easy to get a solution of (\ref{specialb}) with $\mathbf{a}=O(\xi_1^{-\alpha})$ for any $0<\alpha<1$.

\qed

\bigskip

\noindent {\bf Acknowledgments:} J. Wei is supported by a NSERC grant from Canada.

\end{document}